# Generalization error for multi-class margin classification

**Xiaotong Shen**

*School of Statistics*
*The University of Minnesota*
*Minneapolis, MN 55455*
*e-mail:* xshen@stat.umn.edu

**Lifeng Wang**

*Department of Biostatistics*
*University of Pennsylvania*
*Philadelphia, PA 19104*
*e-mail:* lifwang@mail.med.upenn.edu

**Abstract:** In this article, we study rates of convergence of the generalization error of multi-class margin classifiers. In particular, we develop an upper bound theory quantifying the generalization error of various large margin classifiers. The theory permits a treatment of general margin losses, convex or nonconvex, in presence or absence of a dominating class. Three main results are established. First, for any fixed margin loss, there may be a trade-off between the ideal and actual generalization performances with respect to the choice of the class of candidate decision functions, which is governed by the trade-off between the approximation and estimation errors. In fact, different margin losses lead to different ideal or actual performances in specific cases. Second, we demonstrate, in a problem of linear learning, that the convergence rate can be arbitrarily fast in the sample size $n$ depending on the joint distribution of the input/output pair. This goes beyond the anticipated rate $O(n^{-1})$. Third, we establish rates of convergence of several margin classifiers in feature selection with the number of candidate variables $p$ allowed to greatly exceed the sample size $n$ but no faster than $\exp(n)$.

## Contents

## 1. Introduction

Large margin classification has seen significant developments in the past several years, including many well-known classifiers such as Support Vector Machine (SVM, (7)) and Neural Networks. For margin classifiers, this article investigates their generalization accuracy in multi-class classification.

In the literature, the generalization accuracy of large margin classifiers has been investigated in two-class classification. Relevant results can be found in, for example, (3), (29) and (14). For multi-class classification, however, there are many distinct generalizations of the same two-class margin classifier; see Section 3 for a further discussion of this aspect. As a result, much less is known with regard to the generalization accuracy of large margin classifiers, particularly its relation to presence/absence of a dominating class, which is not of concern in the two-class case. Consistency has been studied in (30), and (21). To our knowledge, rates of convergence of the generalization error have not been yet studied for general margin classifiers in multi-class classification.

In the two-class case, the generalization accuracy of a large margin classifier is studied through the notion of Fisher consistency (cf., (15); (30)), where the Bayesian regret $Regret(\hat{f}, \bar{f})$ is used to measure the discrepancy between an estimated decision function $\hat{f}$ and $\bar{f}$, the (global) Bayes decision function over all possible candidate functions. When a specific class of candidate decision functions $\boldsymbol{\mathcal{F}}$ and a surrogate loss $V$ are used in classification, $\bar{f}$ is often not the risk minimizer defined by $V$ over $\boldsymbol{\mathcal{F}}$. Then an approximation error of $\boldsymbol{\mathcal{F}}$ to $\bar{f}$ with respect to $V$ is usually assumed to yield an upper bound of $Regret(\hat{f}, \bar{f})$, expressed in terms of an approximation error plus an estimation error of estimating the decision function. One major difficulty with this formulation is that the approximation error may dominate the corresponding estimation error and be non-zero. This occurs in classification with linear decision functions; see Section 5.1 for an example. In such a situation, well-established bounds for the estimation error become irrelevant, and hence that such a learning theory breaks down when the approximation error does not tend to zero.

To treat the multi-class margin classification, and circumvent the aforementioned difficulty, we take a novel approach by targeting at $Regret(\hat{f}, f^V)$ with $f^V$ the risk minimizer over $\boldsymbol{\mathcal{F}}$ given $V$. Toward this end, we study the ideal generalization performance of $f^V$ and the mean-variance relationship of the cost function. This permits a comparison of various margin classifiers with respect to the ideal and actual performances respectively described in Sections 3 and 4, bypassing the requirement of studying the Fisher consistency. As illustrated in Section 5.2, we show that the rate of convergence of the generalization error of certain large margin classifiers can be arbitrarily fast in linear classification, depending on the joint distribution of the input/output pair. Moreover, in linear classification, the ideal generalization performance is more crucial than the actual generalization performance, whereas in nonlinear classification the approximation error becomes important to the actual generalization performance. Finally, we treat variable selection in sparse learning in a high-dimensional situation. There the focus has been on how to utilize the sparseness structure to attack the curse of high dimensionality, c.f, (31) and (12). Our formulation permits the number of candidate variables $p$ greatly exceeding the sample size $n$. Specifically, we obtain results for several margin classifiers involving feature selection, when $p$ grows no faster than $\exp(n)$. This illustrates the important role of penalty in sparse learning.

This article is organized as follows. Section 2 introduces the notation of generalized multi-class margin losses to unify various generalizations of two-class margin losses. Section 3 discusses the ideal generalization performance of $f^V$ with respect to $V$, whereas Section 4 establishes an upper bound theory concerning the generalization error for margin classifiers. Section 5 illustrates the general theory through four classification examples. The Appendix contains technical proofs.

## 2. Multi-class and generalized margin losses

In $k$-class classification, a decision function vector $\boldsymbol{f} = (f_1, \cdots, f_k)$, with $f_j$ representing class $j$, mapping from input space $\mathcal{X} \subset \mathbb{R}^d$ to $\mathbb{R}$, is estimated through a training sample $Z_i = (\boldsymbol{X}_i, Y_i)_{i=1}^n$, independent and identically distributed according to an unknown joint probability $P(\boldsymbol{x}, y)$, where $Y_i$ is coded as $\{1, \cdots, k\}$. For an instance $\boldsymbol{x}$, classification is performed by rule $\arg\max_{1\le j\le k} f_j(\boldsymbol{x})$, assigning $\boldsymbol{x}$ to a class with the highest value of $f_j(\boldsymbol{x})$; $j = 1, \cdots, k$. The classifier defined by $\arg\max_{1\le j\le k} f_j(\boldsymbol{x})$ partitions $\mathcal{X}$ into $k$ disjoint and exhaustive regions $\mathcal{X}_1, \cdots, \mathcal{X}_k$. To avoid redundancy in $\boldsymbol{f}$, a zero-sum constraint $\sum_{j=1}^k f_j = 0$ is enforced. Note that $f_j$; $j = 1, \cdots, k$ are not probabilities.

In multi-class margin classification, there are a number of generalizations of the same two-class method. We now introduce a framework using the notion of generalized margin, unifying various generalizations. Define the generalized functional margin $\boldsymbol{u}(\boldsymbol{f}(\boldsymbol{x}), y)$ as $(f_y(\boldsymbol{x}) - f_1(\boldsymbol{x}), \ldots, f_y(\boldsymbol{x}) - f_{y-1}(\boldsymbol{x}), f_y(\boldsymbol{x}) - f_{y+1}(\boldsymbol{x}), \ldots, f_y(\boldsymbol{x}) - f_k(\boldsymbol{x})) \equiv (u_1, \cdots, u_{k-1})$, comparing class $y$ against the remaining classes. When $k = 2$, it reduces to the binary functional margin $f_y - f_{c\neq y}$, which, together with the zero-sum constraint, is equivalent to $yf(\boldsymbol{x})$ with $y = \pm 1$. Within this framework, we define a generalized margin loss

$$V(\boldsymbol{f}, z) = h(\boldsymbol{u}(\boldsymbol{f}(\boldsymbol{x}), y))$$

for some measurable function $h$ and $z = (\boldsymbol{x}, y)$, where $V$ is called large margin if it is nondecreasing with respect to each component of $\boldsymbol{u}(\boldsymbol{f}(x), y)$, and $V$ is often called a surrogate loss when it is not the 0-1 loss. For import vector machine (33),

$$h(\boldsymbol{u}) \equiv h_{logit}(\boldsymbol{u}) = \log(1 + \sum_{j=1}^{k-1} \exp(-u_j)).$$

For multi-class SVMs, several versions of generalized hinge loss exist. The generalized hinge loss proposed by (23), (26), (5), (8), and (11) is defined by

$$h(\boldsymbol{u}) \equiv h_1(\boldsymbol{u}) = \sum_{j=1}^{k-1} [1 - u_j]_+;$$

the generalized hinge loss in (13) is defined by

$$h(\boldsymbol{u}) \equiv h_{svm2}(\boldsymbol{u}) = \sum_{j=1}^{k-1} [\frac{\sum_{c=1}^{k-1} u_c}{k} - u_j + 1]_+;$$

the generalized hinge loss in (16) is defined by

$$h(\boldsymbol{u}) \equiv h_{svm3}(\boldsymbol{u}) = [1 - \min_{\{1 \le j \le k-1\}} u_j]_+.$$

For multi-class $\psi$-learning, $h(\boldsymbol{u}) \equiv h_\psi(\boldsymbol{u}) = \psi(\min_{\{1 \le j \le k-1\}} u_j)$ is the generalized $\psi$-loss by (16), with $\psi(x) = 0$ for $x > 1$; 2 for $x < 0$; $2(1-x)$ for $0 \le x \le 1$. For multi-class boosting (32), $h(\boldsymbol{u}) \equiv h_{l_2}(\boldsymbol{u}) = (1 - \min_{\{1 \le j \le k-1\}} u_j)^2$ is a generalized squared loss. Interestingly, for the 0-1 loss, $h(\boldsymbol{u}) \equiv L(\boldsymbol{f}, Z) = I[\min_{\{1 \le j \le k-1\}} u_j < 0]$.

For classification, a penalized cost function is constructed through $V(\boldsymbol{f}, Z)$:

$$n^{-1} \sum_{i=1}^{n} V(\boldsymbol{f}, Z_i) + \lambda J(\boldsymbol{f}) \tag{2.1}$$

where $J(\boldsymbol{f})$ is a nonnegative penalty penalizing undesirable properties of $\boldsymbol{f}$, and $\lambda > 0$ is a tuning parameter controlling the trade-off between training and $J(\boldsymbol{f})$. The minimizer of (2.1) with respect to $\boldsymbol{f} \in \boldsymbol{\mathcal{F}} = \{(f_1, \cdots, f_k) \in \boldsymbol{\mathcal{F}} : \sum_{j=1}^{k} f_j = 0\}$, a class of candidate decision function vectors, yields $\hat{\boldsymbol{f}} = (\hat{f}_1, \cdots, \hat{f}_k)$ thus classifier $\arg\max_{j=1,\cdots,k} \hat{f}_j$.

In classification, $J(\boldsymbol{f})$ is often the inverse of the geometric margin defined by various norms or the conditional Fisher information (6). For instance, in linear SVM classification with feature selection, the inverse geometric margin with respect to a linear decision function vector $\boldsymbol{f}$ is defined as $\frac{1}{2} \sum_{j=1}^{k} \|\boldsymbol{w}_j\|_1$, cf., (4), where $f_j(\boldsymbol{x}) = \langle \boldsymbol{w}_j, \boldsymbol{x} \rangle + b_j$, $j = 1, \cdots, k$, with $\langle \cdot, \cdot \rangle$ the usual inner product in $\mathbb{R}^d$, $b_j \in \mathbb{R}$, and $\| \cdot \|_1$ is the usual $L_1$ norm. In standard kernel SVM learning, the inverse geometric margin becomes $\frac{1}{2} \sum_{j=1}^{k} \|g_j\|_K^2 = \frac{1}{2} \sum_{i=1}^{n} \sum_{k=1}^{n} \alpha_i^j \alpha_k^j K(\boldsymbol{x}_i, \boldsymbol{x}_j)$, where $f_j$ has a kernel representation of $g_j(\boldsymbol{x}) + b_j \equiv \sum_{i=1}^{n} \alpha_i^j K(\boldsymbol{x}, \boldsymbol{x}_i) + b_j$. Here $K(\cdot, \cdot)$ is symmetric and positive semi-definite, mapping from $\mathcal{X} \times \mathcal{X}$ to $\mathbb{R}$, and is assumed to satisfy Mercer's condition (17) so that $\|g\|_K^2$ is a norm.

## 3. Ideal generalization performance

The *generalization error* (GE) is often used to measure the generalization accuracy of a classifier defined by $\boldsymbol{f}$, which is

$$Err(\boldsymbol{f}) = P(Y \neq \arg \max_{j=1,\cdots,k} f_j(\boldsymbol{X})) = EL(\boldsymbol{f}, Z),$$

with multi-class misclassification (0-1) loss $L(\boldsymbol{f}, z) = I(Y \neq \arg \max_{j=1,\cdots,k} f_j(\boldsymbol{X}))$. The corresponding *empirical generalization error* (EGE) is $n^{-1} \sum_{i=1}^{n} L(\boldsymbol{f}, Z_i)$.

Often a surrogate loss $V$ is used in (2.1) as opposed to the 0-1 loss for a computational consideration. In such a situation, (2.1) targets at the minimizer $\boldsymbol{f}^V = \arg \inf_{f \in \boldsymbol{\mathcal{F}}} EV(\boldsymbol{f}, Z)$, which may not belong to $\boldsymbol{\mathcal{F}}$. Consequently, $EV(\boldsymbol{f}^V, Z)$ represents the ideal performance under $V$, whereas $EL(\boldsymbol{f}^V, Z)$ is the ideal generalization performance of $f^V$ when $V$ is used in (2.1). Now define, for $f \in \boldsymbol{\mathcal{F}}$,

$$\begin{aligned} e(\boldsymbol{f}, \boldsymbol{f}^V) &= EL(\boldsymbol{f}, Z) - EL(\boldsymbol{f}^V, Z), \\ e_V(\boldsymbol{f}, \boldsymbol{f}^V) &= EV(\boldsymbol{f}, Z) - EV(\boldsymbol{f}^V, Z). \end{aligned}$$

Note that for $f \in \boldsymbol{\mathcal{F}}$, $e_V(\boldsymbol{f}, \boldsymbol{f}^V) \geq 0$ but $e(\boldsymbol{f}, \boldsymbol{f}^V)$ may not be so, depending on the choice of $V$. In this article, we provide a bound of $|e(\boldsymbol{f}, \boldsymbol{f}^V)|$ to measure the discrepancy between the actual performance and ideal performance of a classifier defined by $\boldsymbol{f}$ in generalization.

It is worthwhile to mention that for two margin losses $V_i$; $i = 1, 2$, the ideal generalization performances determine the asymptotic behavior of their actual generalization performances of the corresponding classifiers defined by $\hat{\boldsymbol{f}}_i$. Therefore, if $EL(\boldsymbol{f}^{V_1}, Z) < EL(\boldsymbol{f}^{V_2}, Z)$ then $EL(\hat{\boldsymbol{f}}_1, Z) < EL(\hat{\boldsymbol{f}}_2, Z)$ eventually provided that $|e(\hat{\boldsymbol{f}}_i, \boldsymbol{f}^{V_i})| \to 0$ as $n \to \infty$. Consequently, a comparison of $|e(\hat{\boldsymbol{f}}_1, \boldsymbol{f}^{V_1})|$ with $|e(\hat{\boldsymbol{f}}_2, \boldsymbol{f}^{V_2})|$ is useful only when their ideal performances are the same, that is, $EL(\boldsymbol{f}^{V_1}, Z) = EL(\boldsymbol{f}^{V_2}, Z)$.

To study the ideal generalization performance of $f^V$ with respect to $V$, let $\bar{\boldsymbol{f}}$ be the (global) Bayes rule, obtained by minimizing $Err(\boldsymbol{f})$ with respect to all $\boldsymbol{f}$, including $\boldsymbol{f} \notin \boldsymbol{\mathcal{F}}$. Note that the (global) Bayes rule is not unique but its error is unique with respect to loss $L$, because any $\bar{\boldsymbol{f}}$, satisfying $\text{argmax}_j \bar{f}_j(x) = \text{argmax}_j P_j(x)$ with $P_j(x) = P(Y = j|X = x)$, yields the same minimal. Without loss of generality, we define $\bar{\boldsymbol{f}} = (\bar{f}_1, \ldots, \bar{f}_k)$ with $\bar{f}_l(x) = \frac{k-1}{k}$ if $l = \text{argmax} P_j(x)$, and $-\frac{1}{k}$ otherwise.

Let $V_{svmj}$ and $V_\psi$ be margin losses defined by $h_{svmj}$ and $h_\psi$, respectively.

**Lemma 1.** *If $\boldsymbol{\mathcal{F}}$ is a linear space, then*

$$EL(\boldsymbol{f}^V, Z) \geq EL(\boldsymbol{f}^{V_\psi}, Z) = EL(\boldsymbol{f}^L, Z) \geq EL(\bar{\boldsymbol{f}}, Z),$$

*for any margin loss $V$. If, in addition, for generalized hinge losses $V_{svmj}$, $j \in \{1, 3\}$, it is separable in that $EV(\boldsymbol{f}^{V_{svmj}}, Z) = 0$, then*

$$EL(\boldsymbol{f}^V, Z) \geq EL(\boldsymbol{f}^{V_{svmj}}, Z) = EL(\boldsymbol{f}^{V_\psi}, Z) = EL(\boldsymbol{f}^L, Z) \geq EL(\bar{\boldsymbol{f}}, Z).$$

Lemma 1 concerns $V_\psi$ in both the separable and nonseparable cases, and $V_{svmj}$; $j = 1, 3$, in the separable case, in relation to other margin losses. For other margin losses, such an inequality may not hold generally, depending on $\boldsymbol{\mathcal{F}}$ and $V$. Therefore a case by case examination may be necessary; see Section 5.1 for an example.

## 4. Actual generalization performance

In our formulation, $\boldsymbol{\mathcal{F}}$ is allowed to depend on the sample size $n$; so is $\boldsymbol{f}^V$ defined by $\mathcal{F}$. When $\boldsymbol{f}^V$ depends on $n$ and approximates $\boldsymbol{f}^*$ (independent of $n$) in that $|e_V(\boldsymbol{f}^V, \boldsymbol{f}^*)| \to 0$ as $n \to \infty$, it seems sensible to use $|e(\boldsymbol{f}, \boldsymbol{f}^*)|$ to measure the actual performance as opposed to $|e(\boldsymbol{f}, \boldsymbol{f}^V)|$. Without loss of generality, we assume that $V \geq 0$.

Let $\boldsymbol{f}_0 = \boldsymbol{f}^*$ when $\boldsymbol{f}^* \in \boldsymbol{\mathcal{F}}$; otherwise $\boldsymbol{f}_0 \in \boldsymbol{\mathcal{F}}$ is chosen such that $e_V(\boldsymbol{f}_0, \boldsymbol{f}^*) \leq \varepsilon_n^2/4$ with $\varepsilon_n$ defined in Assumption C. Now define truncated $V$-loss $V^T$ to be

$$V^T(\boldsymbol{f}, Z) = V(\boldsymbol{f}, Z) \wedge T,$$

for any $\boldsymbol{f} \in \boldsymbol{\mathcal{F}}$ and some truncation constant $T > 0$, where $\wedge$ defines the minimum. Define

$$e_{V^T}(\boldsymbol{f}, \boldsymbol{f}^*) = E(V^T(\boldsymbol{f}, \boldsymbol{Z}) - V^T(\boldsymbol{f}^*, Z)).$$

The following conditions are assumed based on the bracketing $L_2$ metric entropy and the uniform entropy.

**Assumption A:** (Conversion) There exists a constant $T > 0$ independent of $n$ such that $T > \max(V(\boldsymbol{f}_0, Z), V(\boldsymbol{f}^*, Z))$ *a.s.*, and there exist constants $0 < \alpha \leq \infty$ and $c_1 > 0$ such that for all $0 < \epsilon \leq T$ and $\boldsymbol{f} \in \boldsymbol{\mathcal{F}}$,

$$\sup_{\{\boldsymbol{f} \in \mathcal{F}: e_{V^T}(\boldsymbol{f}, \boldsymbol{f}^*) \leq \epsilon\}} |e(\boldsymbol{f}, \boldsymbol{f}^*)| \leq c_1 \epsilon^\alpha.$$

**Assumption B:** (Variance) For some constant $T > 0$, there exist constants $\beta \geq 0$ and $c_2 > 0$ such that for all $0 < \epsilon \leq T$ and $\boldsymbol{f} \in \boldsymbol{\mathcal{F}}$,

$$\sup_{\{\boldsymbol{f} \in \boldsymbol{\mathcal{F}}: e_{V^T}(\boldsymbol{f}, \boldsymbol{f}^*) \leq \epsilon\}} Var(V^T(\boldsymbol{f}, Z) - V(\boldsymbol{f}^*, Z)) \leq c_2 \epsilon^\beta.$$

To specify Assumption C, we define the $L_2$-bracketing metric entropy and the uniform metric entropy for a function space $\boldsymbol{\mathcal{G}} = \{g\}$ consisting of function $g$'s. For any $\epsilon > 0$, call $\{(g_1^l, g_1^u), \ldots, (g_m^l, g_m^u)\}$ an $\epsilon$-bracketing set of $\boldsymbol{\mathcal{G}}$ if for any $g \in \boldsymbol{\mathcal{G}}$ there exists an $j$ such that $g_j^l \leq g \leq g_j^u$ and $\|g_j^u - g_j^l\|_2 \leq \epsilon$, where $\|g\|_2 = (Eg^2)^{1/2}$ is the usual $L_2$-norm. The metric entropy $H_B(\epsilon, \boldsymbol{\mathcal{G}})$ of $\boldsymbol{\mathcal{G}}$ with bracketing is then defined as the logarithm of the cardinality of $\epsilon$-bracketing set of $\boldsymbol{\mathcal{G}}$ of the smallest size. Similarly, a set $(g_1, \cdots, g_m)$ is called an $\varepsilon$-net of $\boldsymbol{\mathcal{G}}$, if for any $g \in \boldsymbol{\mathcal{G}}$, there exists an $j$ such that $\|g_j - g\|_{Q,2} \leq \varepsilon$, where $\|\cdot\|_{Q,2}$ is the $L_2(Q)$-norm with respect to $Q$, defined as $\|g\|_{Q,2} = (\int g^2 dQ)^{1/2}$. The $L_2(Q)$-metric entropy $H_Q(\varepsilon, \mathcal{G})$ is the logarithm of the covering number—minimal size of all $\varepsilon$-nets. The uniform metric entropy is defined as $H_U(\varepsilon, \boldsymbol{\mathcal{G}}) = \sup_Q H_Q(\varepsilon, \boldsymbol{\mathcal{G}})$.

Let $J_0 = \max(J(\boldsymbol{f}_0), 1)$, and $\mathcal{F}_V(s) = \{Z \mapsto V^T(\boldsymbol{f}, Z) - V(\boldsymbol{f}_0, Z) : \boldsymbol{f} \in \boldsymbol{\mathcal{F}}, J(\boldsymbol{f}) \leq J_0 s\}$.

**Assumption C:** (Complexity) For some constants $c_i > 0$; $i = 3, \cdots, 5$, there exists $\varepsilon_n > 0$ such that

$$\sup_{\{s \geq 1\}} \phi(\varepsilon_n, s) \leq c_3 n^{1/2}, \tag{4.1}$$

where $\phi(\varepsilon_n, s) = \int_{c_5 L}^{c_4^{1/2} L^{\beta/2}} H^{1/2}(u, \mathcal{F}_V(s)) du/L$ and $L = L(\varepsilon_n, \lambda, s) = \min(\varepsilon_n^2 + \lambda J_0(s/2 - 1), 1)$, where $H(\cdot, \cdot)$ is $H_B(\cdot, \cdot)$ or $H_U(\cdot, \cdot)$.

Assumption A specifies a relationship between $e(\boldsymbol{f}, \boldsymbol{f}^*)$ and $e_{V^T}(\boldsymbol{f}, \boldsymbol{f}^*)$, which is a first moment condition. Assumption B, on the other hand, relates $e(\boldsymbol{f}, \boldsymbol{f}^*)$ to variance of $(V^T(\boldsymbol{f}, Z) - V(\boldsymbol{f}^*, Z))$. Evidently $Var(V^T(\boldsymbol{f}, Z) - V(\boldsymbol{f}^*, Z)) \leq \min(Var(V(\boldsymbol{f}, Z) - V(\boldsymbol{f}^*, Z)), T^2)$, which implies that $\beta = 0$ in the worst case. Exponents $\alpha$ and $\beta$ in Assumptions A and B are critical to determine the speed of convergence of $e(\boldsymbol{f}, \boldsymbol{f}^*)$, although $e_{V^T}(\boldsymbol{f}, \boldsymbol{f}^*)$ may not converge fast. As illustrated in Section 5.2, an arbitrarily fast rate is achievable in large margin linear classification, because $\alpha$ can be arbitrarily large. Assumption B appears to be important in discriminating several classifiers in the linear and non-linear cases. Assumption C measures the complexity of $\boldsymbol{\mathcal{F}}$. However, if $c_1$ and $c_2$ in Assumptions B and C depend on $n$, then they may enter into the rate.

Two situations are worthwhile mentioning, depending on richness of $\boldsymbol{\mathcal{F}}$. First, when $\boldsymbol{\mathcal{F}}$ is rich, $\boldsymbol{f}^* = \bar{\boldsymbol{f}}$, and margin classification depends only on the behavior of the marginal distribution of $X$ near the decision boundary. This is characterized by the values of $\alpha$ and $\beta$. For instance, in nonlinear multi-class $\psi$-learning, $\alpha = 1$ and $0 < \beta \leq 1$, cf., (16). This corresponds to the case of the $n^{-1}$ rate in the separable case and $n^{-1/2}$ in the non-separable case, as described in (2). Second, when $\boldsymbol{\mathcal{F}}$ is not rich, as in linear classification, $\boldsymbol{f}^* \neq \bar{\boldsymbol{f}}$ is typically the case, where $\alpha$ and $\beta$ depend heavily on the distribution of $(X, Y)$; see Section 5.2 for an example. As a result, actual generalization performances of various margin classifiers are dominated by different ideal generalization performances; see Section 5.1 for an example.

**Theorem 1.** *If Assumptions A-C hold, then, for any estimated decision function vector $\hat{\boldsymbol{f}}$ defined in (2.1), there exists a constant $c_6 > 0$ depending on $c_1$-$c_5$ such that*

$$P\left(e(\hat{\boldsymbol{f}}, \boldsymbol{f}^*) \geq c_1 \delta_n^{2\alpha}\right) \leq c_7 \exp(-c_6 n (\lambda J_0)^{2-\min(\beta,1)}),$$

*provided that $\lambda^{-1} \geq 2\delta_n^{-2} J_0$, where $c_7 = 3.5$ for the bracketing entropy $H_B(\cdot, \cdot)$ and $c_7 = (1 + (20(1 - \frac{1}{32 c_6 n (\lambda J_0)^{2-\min(\beta,1)}})^{-1})^{1/2})$ for the uniform entropy $H_U(\cdot, \cdot)$, and $\delta_n^2 = \min(\varepsilon_n^2 + 2e_V(\boldsymbol{f}_0, \boldsymbol{f}^*), 1)$.*

**Corollary 1.** *Under the assumptions of Theorem 1, $|e(\hat{\boldsymbol{f}}, \boldsymbol{f}^*)| = O_p(\delta_n^{2\alpha})$, $E|e(\hat{\boldsymbol{f}}, \boldsymbol{f}^*)| = O(\delta_n^{2\alpha})$, provided that $n(\lambda J_0)^{2-\min(\beta,1)}$ is bounded away from zero.*

The rate $\delta_n^2$ is governed by two factors: (1) $\varepsilon_n^2$ determined by the complexity of $\boldsymbol{\mathcal{F}}$ and (2) the approximation error $e_V(\boldsymbol{f}^*, \boldsymbol{f}_0)$ defined by $V$. When $e(\boldsymbol{f}^*, \boldsymbol{f}_0) \neq 0$, there is usually a trade-off between the approximation error and the complexity of $\boldsymbol{\mathcal{F}}$ with respect to the choice of $\boldsymbol{f}_0$; see Section 5.3.

**Remark 1:** The results in Theorem 1 and Corollary 1 continue to hold if the "global" entropy is replaced by its corresponding "local" version; see e.g., (24). That is, $\mathcal{F}_V(s)$ is replaced by $\mathcal{F}_V^1(s) = \mathcal{F}_V(s) \cap \{V^T(\boldsymbol{f}, Z) - V(\boldsymbol{f}_0, Z) : e_V(\boldsymbol{f}, \boldsymbol{f}^*) \leq 2s\}$. The proof requires only a slight modification. The local entropy allows us to avoid to loss of $\log(n)$ in linear classification, although it may not be useful for nonlinear classification.

**Remark 2:** For $\psi$-learning, Theorem 1 may be strengthened by replacing $\mathcal{F}$ by the corresponding set entropy if the problem structure is used; cf., (19).

**Remark 3:** The preceding formulation can be easily extended to the situation of multiple regularizers by replacing $\lambda J(\boldsymbol{f})$ by its vector version, i,e, $\boldsymbol{\lambda}^T \boldsymbol{J}(\boldsymbol{f}) = \sum_{j=1}^{l} \lambda_j J_j(\boldsymbol{f})$ with $\boldsymbol{\lambda} = (\lambda_1, \cdots, \lambda_l)^T$ and $\boldsymbol{J}(\boldsymbol{f}) = (J_1(\boldsymbol{f}), \cdots, J_l(\boldsymbol{f}))$.

## 5. Examples

### *5.1. Linear classification: Ideal and actual performances*

This section illustrates that the ideal generalization performances of various margin classifiers, defined by $\boldsymbol{f}^V$, may differ, dominating the corresponding actual ones, where $e(\boldsymbol{f}^V, \bar{\boldsymbol{f}}) \neq 0$. This reinforces our discussion in Section 3.

Consider, for simplicity, a two-class case with $X$ generated from probability density $q(x) = 2^{-1}(\gamma + 1)|x|^\gamma$ for $x \in [-1, 1]$ and some $\gamma \geq 0$. Given $X = x$, $Y$ is sampled from $\{0, 1\}$ according to $P(Y = 1|X = x)$ that is $\theta_1$ if $x > 0$ and $\theta_2$ otherwise, for constants $\theta_1 > 1/2$, $\theta_2 < 1/2$, and $\theta_1 + \theta_2 \neq 1$. Here decision function vector $\boldsymbol{f}$ is $(f, -f)$ with $\boldsymbol{\mathcal{F}} = \{f = ax + b\}$, $\boldsymbol{u}(\boldsymbol{f}(\boldsymbol{x}), \boldsymbol{y})$ is equivalent to $yf(x)$ for coding $y = \pm 1$ with $J(\boldsymbol{f}) = |a|$.

Four margin losses are compared with respect to their ideal and actual performances measured by $e(\boldsymbol{f}^{V_j}, \bar{\boldsymbol{f}}) \geq 0$ and $|e(\hat{\boldsymbol{f}}, \boldsymbol{f}^{V_j})|$. They are exponential, logistic, hinge and $\psi$ losses, denoted as $V_1 = \exp(-yf(x))$, $V_2 = \log(1 + \exp(-yf(x)))$, $V_3 = [1 - yf(x)]_+$, and $V_4 = I[yf(x) \leq 0] + (1 - yf(x))I[0 < yf(x) \leq 1]$.

To obtain an expression of $e_V(\boldsymbol{f}, \bar{\boldsymbol{f}})$ and $e(\boldsymbol{f}, \bar{\boldsymbol{f}})$, let $R_{V_j}(a, b) = EV_j(\boldsymbol{f}, Z)$ and $R(a, b) = EL(\boldsymbol{f}, Z)$. Let $(a_j^*, b_j^*) = \arg\inf R_{V_j}(a, b)$, $j = 1, \cdots, 4$, and $(\bar{a}, \bar{b}) = \arg\inf R_(a, b)$. The expression of $R_{V_j}(a, b)$ is given in the proof of Lemma 2 of the Appendix, with its properties stated in Lemma 2.

**Lemma 2.** *The minimizer* $(\bar{a}, \bar{b}) = (1, 0)$, $(a_j^*, b_j^*)$, $j = 1, \cdots, 3$ *are finite, and* $(a_4^*, b_4^*) = (\infty, 0)$, *or equivalently,* $R_{V_4}(a, b)$ *attains its minimal as* $a \to +\infty$ *and* $b = 0$.

Based on Lemma 2, we compare the ideal performances $e(\boldsymbol{f}^{V_j}, \bar{\boldsymbol{f}}) = R(a_j^*, b_j^*) - R(1, 0)$ for $j = 1, \cdots, 4$. Since $e(\boldsymbol{f}^{V_j}, \bar{\boldsymbol{f}})$ is not analytically tractable, we provide a numerical comparison in the case of $\theta_1 = 3/4$, $\gamma = 0$ and $\theta_2 \in [1/8, 3/8]$.

As displayed in Figure 1, $e(\boldsymbol{f}^{V_j}, \bar{\boldsymbol{f}})$ decreases as $j$ increases from 1 to 4 with $e(\boldsymbol{f}^{V_4}, \bar{\boldsymbol{f}}) = 0$, indicating that $V_4$ dominates $V_1 - V_3$. Note that $e(\boldsymbol{f}^{V_j}, \bar{\boldsymbol{f}}) = 0$; $j = 1, \cdots, 4$, when $\theta_1 = 3/4$ and $\theta_2 = 1/4$, because of symmetry of $q(x)$ in $x$.

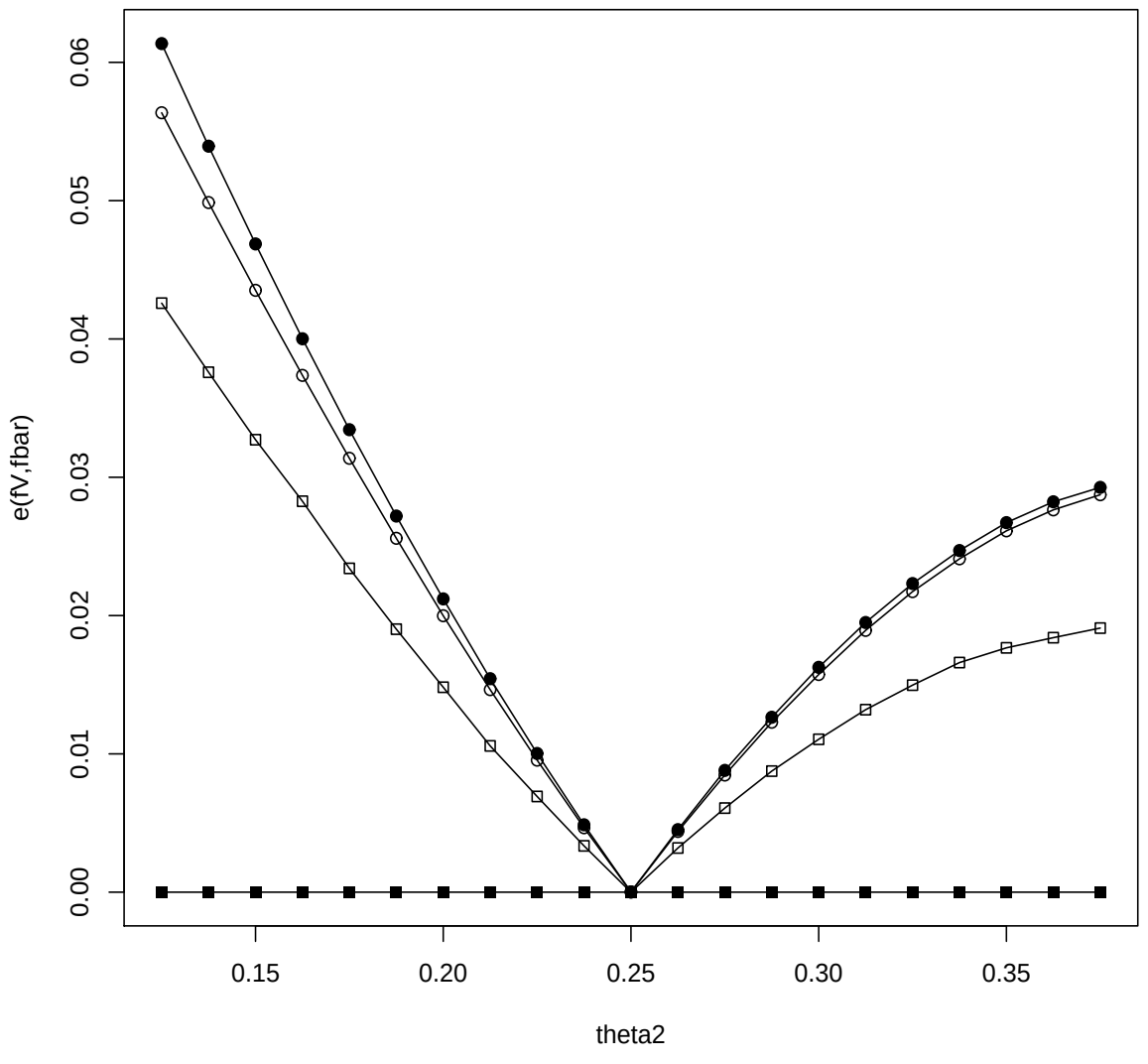

Fig 1. *Plot of $e(\boldsymbol{f}^{V_j}, \bar{\boldsymbol{f}})$ as a function of $\theta_2$, when $\theta_1 = 3/4$ and $\gamma = 0$. Solid circle, circle, square and solid square (from the top curve to the bottom curve) represent $V_j$; $j = 1, \cdots, 4$.*

We now verify Assumptions A-C for $V_j$; $j = 1, \cdots, 4$, with Assumptions A and B checked in Lemma 3.

**Lemma 3.** *Assumptions A and B are for $V_j$, $j = 1, \cdots, 3$, with $\alpha = 1/2$ and $\beta = 1$. For $V_4$, Assumptions A and B hold with $\alpha = 1$ and $\beta = 1$.* To verify Assumption C, let $\boldsymbol{f}_0 = \boldsymbol{f}^{V_j}$ for $V_j$; $j = 1, \cdots, 3$, and compute the local entropy of $\mathcal{F}^1_{V_j}(s) = \mathcal{F}_{V_j}(s) \cap \{V_j(\boldsymbol{f}, \cdot) - V_j(\boldsymbol{f}^{V_j}, \cdot) : e_{V_j}(\boldsymbol{f}, \boldsymbol{f}^{V_j}) \leq 2s\}$. Note that $e_{V_j}(\boldsymbol{f}, \boldsymbol{f}^{V_j}) \leq 2s$ implies that $\|(a,b) - (a^*, b^*)\| \leq c' s^{1/2}$ for the Euclidean norm $\|\cdot\|$ and some $c' > 0$. In addition, for any $g_1, g_2 \in \mathcal{F}^1_{V_j}(s)$, $|g_1(z) - g_2(z)| \leq |(a_1 - a_2)x + (b_1 - b_2)| \leq 2\max(|a_1 - a_2|, |b_1 - b_2|)$. Direct calculation yields that $H_B(u, \mathcal{F}^1_V(s)) \leq c(\log(\min(s^{1/2}, c'u^{1/2})/u^{1/2}))$ for some constant $c > 0$. Easily, $\sup_{s \geq 1} \phi(\varepsilon_n, s) \leq c_4/\varepsilon_n$, when $\lambda \sim \varepsilon_n^2 \leq 1$. Solving (4.1) yields $\varepsilon_n = n^{-1/2}$. By Corollary 1, $e(\hat{\boldsymbol{f}}, \boldsymbol{f}^{V_j}) = O_p(\varepsilon_n^{2\alpha}) = O_p(n^{-1/2})$, and $Ee(\hat{\boldsymbol{f}}, \boldsymbol{f}^{V_j}) = O(n^{-1/2})$ when $\lambda \sim n^{-1}$.

For $V_4$, let $\boldsymbol{f}_0 = (nx, -nx)$. Similarly, $e(\hat{\boldsymbol{f}}, \bar{\boldsymbol{f}}) = O_p(n^{-1})$ and $Ee(\hat{\boldsymbol{f}}, \bar{\boldsymbol{f}}) = O(n^{-1})$ when $\lambda \sim n^{-1}$ for $\psi$-learning.

In conclusion, the ideal performances for $V_1 - V_4$ are usually not equal, with $V_4$ the best as suggested by Lemma 1. The actual performances are dominated by their ideal performances when $\theta_1 + \theta_2 \neq 1$, although $e(\hat{\boldsymbol{f}}, \boldsymbol{f}^{V_j}) = O_p(n^{-1/2})$; $j = 1, \cdots, 3$, and $Ee(\hat{\boldsymbol{f}}, \boldsymbol{f}^{V_4}) = O(n^{-1})$.

### 5.2. Multi-class linear classification: Arbitrarily fast rates

This section illustrates that the rates of convergence for the hinge and logistic losses can be arbitrarily fast even in linear classification. This is because the conversion exponent $\alpha$ in Assumption A can be arbitrarily large, although $e_V(\hat{\boldsymbol{f}}, \boldsymbol{f}^V) = O_P(n^{-1/2})$.

Consider four-class linear classification in which $\boldsymbol{X}$ is sampled according to probability density $q(x_1, x_2) = \lambda \min(|x_1|, |x_2|)^\gamma$ for $(x_1, x_2) \in [-1, 1]^2$, with $\gamma \geq 0$ and normalizing constant $\lambda > 0$. Let $S_j$; $j = 1, \cdots, 4$, be four regions $\{x_1 \geq 0, x_2 \geq 0\}$, $\{x_1 \geq 0, x_2 < 0\}$, $\{x_1 < 0, x_2 \geq 0\}$, and $\{x_1 < 0, x_2 < 0\}$. Now $Y$ is assigned to class $c = j$ with probability $\theta$ $(1/4 < \theta < 1)$ and to the remaining three classes with probability $(1-\theta)/3$ for each, when $x \in S_j$; $j = 1, \cdots, 4$. The (global) Bayes rule $\bar{\boldsymbol{f}} = (x_1 + x_2, -x_1 + x_2, -x_1 - x_2, x_1 - x_2)^T$. In this case, decision function vector $\boldsymbol{f}(\boldsymbol{x})$ is parameterized as $(\boldsymbol{w}_1^T\boldsymbol{x}, \boldsymbol{w}_2^T\boldsymbol{x}, \boldsymbol{w}_3^T\boldsymbol{x}, -(\boldsymbol{w}_1^T + \boldsymbol{w}_2^T + \boldsymbol{w}_3^T)\boldsymbol{x})^T$, defined by $\boldsymbol{w} = (w_{11}, w_{12}, w_{21}, w_{22}, w_{31}, w_{32})$. Here $\boldsymbol{\mathcal{F}}$ consists of such $\boldsymbol{f}$'s and $J(\boldsymbol{f}) = \sum_{c=1}^4 \sum_{j=1}^4 w_{cj}^2$.

For the hinge loss $V(\boldsymbol{f}, z) = h(\boldsymbol{u}(\boldsymbol{f}(\boldsymbol{x}), y))$ with

$$h(\boldsymbol{u}) \equiv h_{svm2}(\boldsymbol{u}) = \sum_{j=1}^{k-1} [\frac{\sum_{c=1}^{k-1} u_c}{k} - u_j + 1]_+,$$

write $EV(\boldsymbol{f}, Z)$ and $EL(\boldsymbol{f}, Z)$ as $R_V(\boldsymbol{w})$ and $R(\boldsymbol{w})$. Then $R_V(\boldsymbol{w})$ is piecewise differentiable, convex, and is minimized by $\boldsymbol{w}^* = r(1, 1, -1, 1, -1, -1)$, where $r$ is the largest negative root of a polynomial in $x$: $9(\theta-1) - 16(\theta-1)x + 12(\theta-1)x^2 + (64\theta - 4)x^4 = 0$. By symmetry of $q(\cdot, \cdot)$, $R_V(\boldsymbol{w})$ is twice-differentiable at $\boldsymbol{w}^*$ with positive definite Hessian matrix $H_1$, implying that $\boldsymbol{f}^V = (f_1^V, \cdots, f_4^V) = ((w_1^*)^T\boldsymbol{x}, (w_2^*)^T\boldsymbol{x}, (w_3^*)^T\boldsymbol{x}, -((w_1^*)^T + (w_2^*)^T + (w_3^*)^T)\boldsymbol{x})^T$, and $e(\boldsymbol{f}^V, \bar{\boldsymbol{f}}) = 0$.

We verify Assumptions A-C, with Assumptions A-B checked in Lemma 4.

**Lemma 4.** *In this example, Assumptions A and B are met for $V$ with $\alpha = 1/2$ and $\beta = 1$.*

For Assumption C, we compute the local entropy of $\mathcal{F}_V^1(s) = \mathcal{F}_V(s) \cap \{V(\boldsymbol{f}, \cdot) - V(\boldsymbol{f}^V, \cdot) : e_V(\boldsymbol{f}, \boldsymbol{f}^V) \leq 2s\}$. Note that $e_V(\boldsymbol{f}, \boldsymbol{f}^V) \leq 2s$ implies that $\|\boldsymbol{w} - \boldsymbol{w}^*\| \leq c's^{1/2}$ for some $c' > 0$, and for any $g, g' \in \mathcal{F}_1(s)$, $|g(\boldsymbol{z}) - g'(\boldsymbol{z})| \leq \sum_{c=1}^4 |f_c(\boldsymbol{x}) - f_c'(\boldsymbol{x})| \leq 12 \max_{1 \leq c \leq 4, 1 \leq j \leq 3}(|w_{cj} - w_{cj}'|)$. Direct calculation yields that $H_B(u, \mathcal{F}_V^1(s)) \leq O(\log(\min(s^{1/2}, c's^{1/2})/u^{1/2}))$, and that $\sup_{s \geq 1} \phi(\varepsilon_n, s) \leq c_4/\varepsilon_n$, when $\lambda \sim \varepsilon_n^2 \leq 1$. Solving (4.1) yields $\varepsilon_n = n^{-1/2}$. By Corollary 1, $e(\hat{\boldsymbol{f}}, \bar{\boldsymbol{f}}) = e(\hat{\boldsymbol{f}}, \boldsymbol{f}^V) = |e(\hat{\boldsymbol{f}}, \boldsymbol{f}^V)| = O_p(n^{-(\gamma+1)/2})$ and $E|e(\hat{\boldsymbol{f}}, \boldsymbol{f}^V)| = n^{-(\gamma+1)/2}$ when $\lambda \sim n^{-1}$.

For the logistic loss, an application of the same argument yields that the same $\boldsymbol{w}^* = r(1, 1, -1, 1, -1, -1)$ as the minimizer of $R_V(\boldsymbol{w})$. Furthermore, $\alpha = (\gamma + 1)/2$ and $\beta = 1$, yielding that the same rates as the hinge loss.

Interestingly, the fast rate $e(\hat{\boldsymbol{f}}, \bar{\boldsymbol{f}}) = e(\hat{\boldsymbol{f}}, \boldsymbol{f}^V) = O_p(n^{-(\gamma+1)/2})$ is because classification is easier than its counterpart-function estimation, as measured by $\gamma \geq 0$. This is evident from that $e_V(\hat{\boldsymbol{f}}, \boldsymbol{f}^V) = O_p(n^{-1/2})$. This rate is arbitrarily fast as $\gamma \to \infty$.

### 5.3. Nonlinear classification: Spline kernels

This section examines one nonlinear learning case and the issue of dominating class with regard to generalization error for multi-class SVM with $V_{svm1} \equiv h_{svm1} \equiv \sum_{j=1}^{k-1}[1 - u_j(\boldsymbol{f}(x), y)]_+$ and multi-class $\psi$-learning with $V_\psi \equiv h_\psi = \psi(\min_{\{1 \leq j \leq k-1\}} u_j)$, as defined in Section 2. Consider three-class classification with spline kernels, where $X$ is generated according to the uniform distribution over $[0,1]$. Given $X = x$, $\boldsymbol{Y} = (Y_1, \cdots, Y_3)$ is sampled from $(P(Y = 1|X = x), P(Y = 2|X = x), P(Y = 3|X = x)) \equiv (p_1(x), p_2(x), p_3(x))$, which is $(5/11, 3/11, 3/11)$ when $x \leq 1/3$, $(3/11, 5/11, 3/11)$ when $1/3 < x \leq 2/3$, and $(3/11, 3/11, 5/11)$ when $x > 2/3$. Evidently, for each $x \in [0,1]$, there does not exist a dominating class because $\max_{1 \leq i \leq 3} p_i(x) = 5/11 < 1/2$.

In this example, $\boldsymbol{\mathcal{F}} = \{(f_1, f_2, f_3) : f_i \in W_m[0,1], \sum_{i=1}^3 f_i = 0\}$ is defined by a Sobolev space $W_m[0,1] = \{f : f^{(m-1)}$is absolutely continuous, $f^{(m)} \in L_2[0,1]\}$ with the degree of smoothness $m$ measured by the $L_2$-norm, generated by spline kernel $K(\cdot,\cdot)$ whose expression can be found in (10) or (23). In what follows, we embed $\{(f_1, f_2, f_3) : f_c = \sum_{i=1}^n \alpha_{ic} K(x_i, x)\}$ with penalty $J(\boldsymbol{f}) = \sum_{c=1}^3 \int_0^1 f_c^{(m)}(u)^2 du$ in (2.1) into $\boldsymbol{\mathcal{F}}$ with penalty $J(\boldsymbol{f}) = \sum_{c=1}^3 \sum_{j=0}^{m-1} f_c^{(j)}(0)^2 + \int_0^1 f_c^{(m)}(u)^2 du$. It follows from the reproducing kernel Hilbert spaces (RKHS) representation theorem (cf., (10)) that minimization of (2.1) over $\boldsymbol{\mathcal{F}}$ is equivalent to that over its subspace $\{(f_1, f_2, f_3) : f_c = \sum_{i=1}^n \alpha_{ic} K(x_i, x)\}$ with $J(\boldsymbol{f}) = \sum_{c=1}^3 \int_0^1 f_c^{(m)}(u)^2 du$.

We now verify Assumptions A-C. Some useful facts are given in Lemmas 5-7.

**Lemma 5.** *(Global Bayes rule $\bar{\boldsymbol{f}}$) In this example,* $\bar{\boldsymbol{f}}(x) = (I[0 \leq x \leq 1/3] - 1/3, I[1/3 < x \leq 2/3] - 1/3, I[2/3 < x \leq 1] - 1/3) = \arg\inf_{\boldsymbol{f}} EV(\boldsymbol{f}, \boldsymbol{Z})$ *for* $V = V_{svm1}, V_\psi$.

**Lemma 6.** *(Assumption A) In this case,* $e_{V^T}(\boldsymbol{f}, \boldsymbol{f}^*) \geq Ce(\boldsymbol{f}, \boldsymbol{f}^*)$ *for* $V = V_{svm1}$, *some constant* $C > 0$, *any* $T \geq 9$ *and any measurable* $\boldsymbol{f} \in \mathbb{R}^3$.

**Lemma 7.** *(Assumption B) In this example,* $E(V^T(\boldsymbol{f}, \boldsymbol{Z}) - V(\boldsymbol{f}^*, \boldsymbol{Z}))^2 \leq CE(V^T(\boldsymbol{f}, \boldsymbol{Z}) - V(\boldsymbol{f}^*, \boldsymbol{Z}))$ *for some constant* $C > 0$, *any* $T \geq 9$ *and any measurable* $\boldsymbol{f} \in \mathbb{R}^3$, *where* $V = V_{svm1}$.

By Lemma 5, $\bar{\boldsymbol{f}}(x) = (I[0 < x \leq 1/3] - 1/3, I[1/3 < x \leq 2/3] - 1/3, I[2/3 \leq x \leq 1] - 1/3)$, which can be approximated by $\boldsymbol{\mathcal{F}}$ with respect to the $L_1$-norm $\|\cdot\|_1$ ($\|f\|_1 = E|f(X)|$). Hence $\boldsymbol{f}^* = \bar{\boldsymbol{f}} = \boldsymbol{f}^V = \arg\inf_{\boldsymbol{f}} EV(\boldsymbol{f}, Z)$ for $V = V_{svm1}, V_\psi$.

**SVM:** Let $\boldsymbol{f}_0 = (f^{(1)}, -f^{(1)} - f^{(3)}, f^{(3)})$ with $f^{(1)} = 2/3 - (1 + \exp(-\tau(x - 1/3)))^{-1}$, $f^{(3)} = (1 + \exp(-\tau(x - 1/3)))^{-1} - 1/3$, with parameter $\tau > 0$ to be specified. Set $T = 4$. Then $T \geq \max(\sup_z V_{svm1}(\boldsymbol{f}_0, z), \sup_z V_{svm1}(\boldsymbol{f}^*, z)) \geq 0$. By Lemmas 6 and 7, Assumptions A and B are met with $\alpha = 1$ and $\beta = 1$. For Assumption C, it can be verified that $e_{V_{svm1}}(\boldsymbol{f}_0, \boldsymbol{f}^*) \le 2\int_0^1 \|\boldsymbol{f}_0(u) - \boldsymbol{f}^*(u)\|_1 du = O(\tau^{-1})$, and $J(\boldsymbol{f}_0) = O(\tau^{2m-1})$. By Proposition 6 of (9), $H_B(u, \mathcal{F}_{V_{svm1}}(s)) = O(((J_0 s)^{1/2}/u)^{1/m})$ with $J_0 = O(\tau^{2m-1})$. Solving (4.1) yields a rate $\varepsilon_n^2 = O((J_0^{\frac{1}{4m}} n^{-1/2})^{\frac{4m}{2m+1}})$ when $J_0\lambda \sim \varepsilon_n^2$. As a result, we have $e_{V_{svm1}}(\hat{\boldsymbol{f}}, \boldsymbol{f}^*) = O_P(\max(\varepsilon_n^2, e_{V_1}(\boldsymbol{f}_0, \boldsymbol{f}^*)) = O_P(\max(\tau^{\frac{2m-1}{2m+1}} n^{-\frac{2m}{2m+1}}, \frac{1}{\tau})) = O_P(n^{-1/2})$, with a choice of $\tau \sim n^{1/2}$ and $\lambda \sim n^{-m}$.

**$\psi$-learning:** Set $T \ge 1$ as $0 < V_4 \le 1$. For Assumption A, $\alpha = 1$ by Theorem 3.1 of (16). For Assumption B, $\beta = 1$ following an argument similar to that in (16). For Assumption C, let $f_0(x) = \tau(4 - 9x, 1, 9x - 5)$ when $m \ge 2$; $\boldsymbol{f}_0 = (f^{(1)}, -f^{(1)} - f^{(3)}, f^{(3)})$ with $f^{(1)} = 2/3 - (1 + \exp(-\tau(x - 1/3)))^{-1}$, $f^{(3)} = (1 + \exp(-\tau(x - 1/3)))^{-1} - 1/3$, when $m = 1$. Then $e_{V_\psi}(\boldsymbol{f}_0, \boldsymbol{f}^*) = O(\tau^{-1})$, $J(\boldsymbol{f}_0) = O(\tau^q)$ with $q = 1$ when $m = 1$ and $q = 2$ when $m \ge 2$, and $H_B(u, \mathcal{F}_{V_\psi}(s)) = O(((J_0 s)^{1/2}/u)^{1/m})$. Solving (4.1) yields a rate $\varepsilon_n^2 = O((J_0^{\frac{1}{4m}} n^{-1/2})^{\frac{4m}{2m+1}})$ when $J_0\lambda \sim \varepsilon_n^2$. By Corollary 1, when $m \ge 2$, $e(\hat{\boldsymbol{f}}, \bar{\boldsymbol{f}}) = O_P(\max(\varepsilon_n^2, e_{V_\psi}(\boldsymbol{f}_0, \bar{\boldsymbol{f}}))) = O_P(\max(\tau^{\frac{2}{2m+1}} n^{-\frac{2m}{2m+1}}, \frac{1}{\tau})) = O_P(n^{-2m/(2m+3)})$; when $m = 1$, $e(\hat{\boldsymbol{f}}, \bar{\boldsymbol{f}}) = O_p(n^{-1/2})$ $\tau \sim n^{1/2}$ and $\lambda \sim n^{-1}$. This yields $Ee(\hat{\boldsymbol{f}}, \bar{\boldsymbol{f}}) = O(n^{-2m/(2m+3)})$ when $m \ge 2$, with $\tau = n^{2m/(2m+3)}$ and $\lambda \sim n^{-6m/(2m+3)}$; $Ee(\hat{\boldsymbol{f}}, \bar{\boldsymbol{f}}) = O(n^{-1/2})$ when $m = 1$ with $\tau \sim n^{1/2}$ and $\lambda \sim n^{-1}$.

Evidently, the approximation error $e_V(\boldsymbol{f}_0, \bar{\boldsymbol{f}})$ and $J_0 = \max(J(\boldsymbol{f}_0), 1)$ play a key role in rates of convergence. With different choices of approximating $f_0$ for $\psi$-loss and the hinge loss, $\psi$-learning and SVM have different error rates with the $\psi$-loss yielding a faster rate when $m \ge 2$ and the same rate when $m = 1$. Moreover, in this example, the dominating class does not seem to be an issue.

### *5.4. Feature selection: High-dimension p but low sample size n*

This section illustrates applicability of general theorem to the high-dimension, low sample size situation. Consider feature selection in classification, where the number of candidate covariates $p$ is allow to greatly exceed the sample size $n$ and to depend on $n$. For the $L_1$ penalty, (22) and (28) obtained the rates of convergence for the binary SVM when $p < n$ and multi-class SVM when $p > n$.

Here we apply the general theory to the elastic-net penalty (see (31)) for binary SVM (27) to obtain a parallel result of (28). We use linear representations in (2.1) as in (27), because of over-specification of non-linear representations. Here decision function vector $\boldsymbol{f}$ is $(f, -f)$ with $f \in \boldsymbol{\mathcal{F}} = \{f(\boldsymbol{x}) = \boldsymbol{w}^T\boldsymbol{x} : \boldsymbol{x} \in [-1, 1]^p\}$, and $J(\boldsymbol{f}) = J_\theta(f) = \theta\|\boldsymbol{w}\|_1 + (1 - \theta)\|\boldsymbol{w}\|_2^2$ is a weighted average of the $L_1$ and $L_2$ norms with a weighting parameter $\theta \in [0, 1]$, cf., (31).

In this example, $(\boldsymbol{X} = (X_1, \cdots, X_p), Y)$ are generated as follows. First, randomly sample $\boldsymbol{X}$ according to the uniform distribution $[-1, 1]^p$. Second, given $\boldsymbol{X} = \boldsymbol{x}$, $Y$ is sampled according to $P(Y = 1|X_1 = x_1)$, which is $\tau > 1/2$ if $x_1 > 0$, and $1 - \tau$ if $x_1 \le 0$. This is a version of Example 5.1 with $\gamma = 0$ and $\theta_1 = \tau$ and $\theta_2 = 1 - \tau$ in a high-dimensional situation. Evidently, $(X_2 \cdots, X_p)$ are redundant variables.

We now verify Assumptions A-C for the hinge loss $V$. Because $X_1$ and $Y$ are independent of $(X_2, \dots, X_p)$, one can verify that the minimal of $EV(\mathbf{f}, \mathbf{Z})$ is that of $EV(f_1(X_1), Y)$ over $\{f_1 : f_1(x) = ax_1 + b\}$, attained by $f_1^* = a^* x_1$ for some $a^* > 0$. For Assumptions A-B, we apply the result in Example 5.1 to obtain $\alpha = 1/2$ and $\beta = 1$.

For Assumption C, we apply Lemma 8 to compute $H_U(\varepsilon, \mathcal{F}_V(s))$.

**Lemma 8.** *For $\mathcal{G}_p(s) = \{f(\boldsymbol{x}) = \boldsymbol{w}^T \boldsymbol{x} : \boldsymbol{w}, \boldsymbol{x} \in [-1,1]^p, \|\boldsymbol{w}\|_1 \le s\}$ and any $\varepsilon > 0$, there exists a constant $c > 0$ such that $H_U(\varepsilon, \mathcal{G}_p(s)) \le cs^2(p \log(1 + \frac{1}{p\varepsilon^2}) + \varepsilon^{-2} \log(p\varepsilon^2 + 1))$.*

Note that $\mathcal{F}(s) \subset \mathcal{G}_p(\theta^{-1} s)$, by Lemma 8, $H_U(\varepsilon, \mathcal{F}_V(s)) = O(p \log(1 + \frac{1}{p\varepsilon^2}) + \varepsilon^{-2} \log(p\varepsilon^2 + 1))$. Set $\lambda \sim \varepsilon_n^2/(2J(f^*))$. To solve $\sup_s \phi(\varepsilon_n, s) \le c_2 n^{1/2}$, note that $\sup_s \phi(\varepsilon_n, s) = \phi(\varepsilon_n, s^*)$ for some finite $s^*$. Then it suffices to solve $\phi(\varepsilon_n, s^*) \le c_2 n^{1/2}$, involving

$$
\begin{aligned}
&\int_{c_4 2^{-1}\varepsilon_n^2}^{c_3(2^{-1}\varepsilon_n^2)^{1/2}} (p \log(1 + \frac{1}{pu^2}) + u^{-2} \log(pu^2 + 1))^{1/2} \, du \\
\le \; & \int_{c_4(2^{-1}\varepsilon_n^2)}^{c_3(2^{-1}\varepsilon_n^2)^{1/2}} (p^{1/2} \log^{1/2}(1 + \frac{1}{pu^2}) + u^{-1} \log^{1/2}(pu^2 + 1)) \, du \\
= \; & \int_{c_4^2 p(2^{-1}\varepsilon_n^2)^2}^{c_3^2 p(2^{-1}\varepsilon_n^2)} t^{1/2} \log^{1/2}(1 + \frac{1}{t}) \, dt + \int_{c_4^2 p(2^{-1}\varepsilon_n^2)^2}^{c_3^2 p(2^{-1}\varepsilon_n^2)} t^{-1} \log^{1/2}(t+1) \, dt \\
\equiv \; & I_1 + I_2.
\end{aligned}
$$

Three cases are examined: First, when $p\varepsilon_n^2 = o(1)$, and $I_1 + I_2 \le 2I_1 = O((p\varepsilon_n^2 \log((p\varepsilon_n^2)^{-1}))^{1/2})$. Solving $\phi(\varepsilon_n, a^*) \le c_2 n^{1/2}$ is equivalent to solving $(p\epsilon_n^2)^{1/2} \log^{1/2}((p\varepsilon_n^2)^{-1}) = O(n^1/2\varepsilon_n^2)$ with respect to $\varepsilon_n^2$, which yields $\varepsilon_n^2 = O((p/n) \log(n/p))$. When $p\varepsilon_n^2 = O(1)$, there exist two constants $0 < B_1, B_2 < \infty$ such that $B_1 < I_1 + I_2 < B_2$, implying $\varepsilon_n^2 = O(n^{-1/2})$. When $p\varepsilon_n^2 \to \infty$, $I_1 \ge I_2$ and $I_1 + I_2 \le 2I_2 = O(\log^{1/2}(p\varepsilon_n^2) \log({\varepsilon_n^2}^{-1}))$. Solving equation $\log^{1/2}(p\varepsilon_n^2) \log({\varepsilon_n^2}^{-1}) = O(n^{1/2}\varepsilon_n^2)$ yields $\varepsilon_n^2 = (n^{-1} \log(t_n p))^{1/2} \log(n)$ when $t_n = (n^{-1} \log p)^{1/2} \log(n/\log p)$.

As a result, the rate is $\varepsilon_n^2 = (\frac{p}{n} \log(\frac{n}{p}))^{1/2}$ when $p \ll n^{1/2}$, $\varepsilon_n^2 = n^{-1/2}$ when $p = O(n^{1/2})$, and $\varepsilon_n = [(n^{-1} \log(t_n p))^{1/2} \log(n)]^{1/2}$ with a choice of $t_n = (n^{-1} \log p)^{1/2} \log(\frac{n}{\log p})$ when $p \gg n^{1/2}$ but $\log p/n = o(1)$. Note that in the last case Assumption B plays no role when $\mathcal{F}$ is too large. By Corollary 1, $e(\hat{\boldsymbol{f}}, \bar{\boldsymbol{f}}) = O_p(\varepsilon_n)$ when $\lambda \sim \varepsilon_n^2$.

## 6. Conclusion

This article develops a statistical learning theory for quantifying the generalization error of large margin classifiers in multi-class classification. In particular, the theory develops upper bounds for a general large margin classifier, which permits a theoretical treatment for the situation of high-dimension but low sample size. Through the theory, several learning examples are studied, where the generalization errors for several large margin classifiers are established. In a linear case, fast rates of convergence are obtained, and in a case of sparse learning, rates are derived for feature selection in which the number of variable greatly exceeds the sample size.

To compare various large margin classifiers with regard to generalization, we may need to develop a lower bound theory. Otherwise, a comparison may be inconclusive although our learning theory provides an upper bound result.

**Acknowledgments.** The author would like to thank the reviewers for helpful comments and suggestions. This research was supported in part by National Science Foundation Grants IIS-0328802 and DMS-0604394.

## Appendix A: Technical proofs

**Proof of Lemma 1:** To prove $EL(\boldsymbol{f}^{V_\psi}, Z) = EL(\boldsymbol{f}^L, Z)$, note that it follows from the definition of $\boldsymbol{f}^L$ that $EL(\boldsymbol{f}^V, Z) \geq EL(\boldsymbol{f}^L, Z)$. Then for any $\varepsilon > 0$ there exists $\boldsymbol{f}_0 \in \boldsymbol{\mathcal{F}}$ such that $EL(\boldsymbol{f}_0, Z) \leq EL(\boldsymbol{f}^L, Z) + \varepsilon$. It follows from linearity of $\boldsymbol{\mathcal{F}}$ that $c\boldsymbol{f}_0 \in \boldsymbol{\mathcal{F}}$ for any constant $c > 0$. The result then follows from the fact that $\lim_{c\to\infty} EV_\psi(c\boldsymbol{f}_0, Z) = EL(\boldsymbol{f}_0, Z)$.

It follows from the fact that $EL(\boldsymbol{f}^V, Z) \geq EL(\boldsymbol{f}^L, Z)$ that $EL(\boldsymbol{f}^V, Z) \geq EL(\boldsymbol{f}^{V_\psi}, Z) = EL(\boldsymbol{f}^L, Z)$.

For $V \equiv h_{svmj}$, in the separable case, the result follows from that $h_{svm1}(z) \geq h_{svm3}(z) = \frac{1}{2}V_\psi(z)$ for $z \geq 0$.

**Proof of Theorem 1:** First we introduce some notations to be used. Let $\tilde{V}(\boldsymbol{f}, Z) = V(\boldsymbol{f}, Z) + \lambda J(\boldsymbol{f})$ and $\tilde{V}^T(\boldsymbol{f}, Z) = V^T(\boldsymbol{f}, Z) + \lambda J(\boldsymbol{f})$. Define the scaled empirical process $E_n(V^T(\boldsymbol{f}, Z) - V^T(\boldsymbol{f}_0, Z))$ as $n^{-1}\sum_{i=1}^n (V^T(\boldsymbol{f}, Z_i) - V^T(\boldsymbol{f}_0, Z_i)$. Let $A_{i,j} = \{f \in \boldsymbol{\mathcal{F}} : 2^{i-1}\delta_n^2 \leq e_{V^T}(\boldsymbol{f}, \boldsymbol{f}^*) < 2^i\delta_n^2, 2^{j-1}\max(J_0, 1) \leq J(\boldsymbol{f}) < 2^j \max(J_0, 1)\}$ and $A_{i,0} = \{f \in \boldsymbol{\mathcal{F}} : 2^{i-1}\delta_n^2 \leq e_{V^T}(\boldsymbol{f}, \boldsymbol{f}^*) < 2^i\delta_n^2, J(\boldsymbol{f}) < \max(J_0, 1)\}$, for $j = 1, 2, \cdots$, and $i = 1, 2, \cdots$.

The treatment here is to use a large deviation inequality in Theorem 3 of (20) for the bracketing entropy and Lemma 9 below for the uniform entropy. Our approach for bounding $P\left(|e(\hat{\boldsymbol{f}}, \boldsymbol{f}^*)| \geq \delta_n^2\right)$ is to bound a sequence of empirical processes induced by the cost function $l$ over $P(A_{ij})$; $i, j = 1, \cdots, n$. Specifically, we apply a large deviation inequality for empirical processes, by controlling the mean and variance defined by $V(\boldsymbol{f}, Z_i)$ and penalty $\lambda$. This yields an inequality for empirical processes and thus for $e(\hat{\boldsymbol{f}}, \boldsymbol{f}^*)$. In what follows, we shall prove the case of the bracketing entropy as that for the uniform entropy is essentially the same.

First we establish a connection between $e(\boldsymbol{f}, \boldsymbol{f}^*)$ and the cost function. By the definition of $\hat{\boldsymbol{f}}$, $\boldsymbol{f}_0$ ($e_{V^T}(\boldsymbol{f}_0, \boldsymbol{f}^*) \leq e_V(\boldsymbol{f}_0, \boldsymbol{f}^*) < \delta_n^2$), and Assumption A, $\left\{|e(\hat{\boldsymbol{f}},\boldsymbol{f}^*)| \geq c_1\delta_n^{2\alpha}\right\} \subset \left\{e_{V^T}(\hat{\boldsymbol{f}},\boldsymbol{f}^*) \geq \delta_n^2\right\}$ is a subset of

$$\begin{aligned} &\left\{ \sup_{\{f\in\boldsymbol{\mathcal{F}}:e_{V^T}(\boldsymbol{f},\boldsymbol{f}^*)\geq\delta_n^2\}} \sum_{i=1}^n (\tilde{V}(\boldsymbol{f}_0,Z_i) - \tilde{V}(\boldsymbol{f},Z_i)) \geq 0 \right\} \\ \subset\ &\left\{ \sup_{\{f\in\boldsymbol{\mathcal{F}}:e_{V^T}(\boldsymbol{f},\boldsymbol{f}^*)\geq\delta_n^2\}} \sum_{i=1}^n (\tilde{V}(\boldsymbol{f}_0,Z_i) - \tilde{V}^T(\boldsymbol{f},Z_i)) \geq 0 \right\}. \end{aligned}$$

Hence $P\left(|e(\hat{\boldsymbol{f}},\boldsymbol{f}^*)| \geq c_1\delta_n^{2\alpha}\right)$ is upper bounded by

$$\begin{aligned} I &\equiv P^*\left( \sup_{\{\boldsymbol{f}\in\boldsymbol{\mathcal{F}}:e_{V^T}(\boldsymbol{f},\boldsymbol{f}^*)\geq\delta_n^2\}} n^{-1}\sum_{i=1}^n (\tilde{V}(\boldsymbol{f}_0,Z_i) - \tilde{V}^T(\boldsymbol{f},Z_i)) \geq 0 \right), \\ &\leq I_1 + I_2, \end{aligned}$$

where $P^*$ is the outer probability, and

$$\begin{aligned} I_1 &= \sum_{i,j\geq 1} P^*\left( \sup_{\boldsymbol{f}\in A_{i,j}} n^{-1}\sum_{i=1}^n (\tilde{V}(\boldsymbol{f}_0,Z_i) - \tilde{V}^T(\boldsymbol{f},Z_i)) \geq 0 \right) \\ I_2 &= \sum_{i=1}^{\infty} P^*\left( \sup_{\boldsymbol{f}\in A_{i,0}} n^{-1}\sum_{i=1}^n (\tilde{V}(\boldsymbol{f}_0,Z_i) - \tilde{V}^T(\boldsymbol{f},Z_i)) \geq 0 \right). \end{aligned}$$

To bound $I_1$, consider $P^*\left(\sup_{\boldsymbol{f}\in A_{i,j}} n^{-1}\sum_{i=1}^n(\tilde{V}(\boldsymbol{f}_0,Z_i) - \tilde{V}^T(\boldsymbol{f},Z_i)) \geq 0\right)$, for each $i = 1,\cdots, j = 0,\cdots$. Let $M(i,j) = 2^{i-1}\delta_n^2 + \lambda 2^{j-1}J(\boldsymbol{f}_0)$. For the mean, using the assumption that $\delta^2/2 \geq \lambda J_0$ and the fact that $e_{V^T}(\boldsymbol{f}_0,\boldsymbol{f}^*) = e_V(\boldsymbol{f}_0,\boldsymbol{f}^*) < \delta_n^2/2$, it follows that $\inf_{A_{i,j}} \mathrm{E}(\tilde{V}^T(\boldsymbol{f},Z_1) - \tilde{V}(\boldsymbol{f}_0,Z_1))$ is lower bounded by

$$\inf_{A_{i,j}} \mathrm{E}(V^T(\boldsymbol{f},Z_1) - V(\boldsymbol{f}^*,Z_1) + \lambda(J(\boldsymbol{f}) - J(\boldsymbol{f}_0))) - e_{V^T}(\boldsymbol{f}_0,\boldsymbol{f}^*) \geq M(i,j),$$

$i = 1,\cdots, j = 0,\cdots$. Similarly, for the variance, it follows from Assumption B and the fact that $Var(V^T(\boldsymbol{f},Z_1) - V(\boldsymbol{f}_0,Z_1)) \leq 2[Var(V^T(\boldsymbol{f},Z_1) - V(\boldsymbol{f}^*,Z_1)) + Var(V(\boldsymbol{f}_0,Z_1) - V(\boldsymbol{f}^*,Z_1))]$ that

$$\sup_{A_{i,j}} Var(V^T(\boldsymbol{f},Z_1) - V(\boldsymbol{f}_0,Z_1)) \leq 4c_2M^\beta(i,j);$$

$i = 1,\cdots, j = 0,\cdots,$.

Note that $0 < \delta_n \leq 1$ and $\lambda\max(J_0,1) \leq \delta_n^2/2$. An application of Theorem 3 of (20) with $M = n^{1/2}M(i,j)$, $v = 4c_2M^\beta(i,j)$, $\varepsilon = 1/2$, and $|V(\boldsymbol{f}_0,Z_i) -$

$V^T(\boldsymbol{f}, Z_i)| \leq 2T$, yields, by Assumption C, that

$$\begin{aligned} I_1 &\leq \sum_{i,j:M(i,j)\leq T} 3\exp\left(-\frac{(1-\varepsilon)nM(i,j)^2}{2(4M^\beta(i,j)+M(i,j)T/3)}\right) \\ &\leq \sum_{j=1}^\infty\sum_{i=1}^\infty 3\exp(-c_6 nM(i,j)^{2-\min(1,\beta)}) \\ &\leq \sum_{j=1}^\infty\sum_{i=1}^\infty 3\exp(-c_6 n[2^{i-1}\delta_n^2+(2^{j-1}-1)\lambda J_0]^{2-\min(1,\beta)}) \\ &\leq \sum_{j=1}^\infty\sum_{i=1}^\infty 3\exp(-c_6 n[(2^{i-1}\delta_n^2)^{2-\min(1,\beta)}+((2^{j-1}-1)\lambda J_0)^{2-\min(1,\beta)}]) \\ &\leq 3\exp(-c_6 n(\lambda J_0)^{2-\min(1,\beta)})/[(1-\exp(-c_6 n(\lambda J_0)^{2-\min(1,\beta)}))]^2. \end{aligned}$$

Here and in the sequel $c_6$ is a positive generic constant. Similarly, $I_2$ can be bounded.

To prove the result with the uniform entropy, we use Lemma 9 with a slight modification of the proof.

Finally,

$$I \leq I_1 + I_2 \leq 6\exp(-c_6 n(\lambda J_0)^{2-\min(1,\beta)})/[(1-\exp(-c_5 n(\lambda J_0)^{2-\min(1,\beta)}))]^2.$$

This implies that $I^{1/2} \leq (5/2+I^{1/2})\exp(-c_6 n(\lambda J_0)^{2-\min(1,\beta)})$. The result then follows from the fact $I \leq I^{1/2} \leq 1$. □

Now we derive Lemma 9 as a version of Theorem 1 of (20) using the uniform entropy.

**Lemma 9.** *Let $\mathcal{F}$ be a collection of functions $f$ with $0 \leq f \leq 1$, $P_n(f) = n^{-1}\sum_{i=1}^n f(Y_i)$, $Pf = Ef(Y_1) = 0$, $Y_i \sim$ i.i.d, and let $v > \sup_{f\in\mathcal{F}} Pf^2 = \sup_{f\in\mathcal{F}} Var(f)$. For $M > 0$ and real $\theta \in (0,1)$, let $\psi(M,n,v) = \frac{nM^2}{128v}$ and $s = \frac{\theta M}{8\sqrt{6}}$. Suppose*

$$H_U(v^{1/2},\mathcal{F}) \leq \frac{\theta}{4}\psi(M,n,v), \tag{A.1}$$

$$M \leq 16(1-3\theta/4)^{1/2}v, \tag{A.2}$$

*and, if $s \leq v^{1/2}$,*

$$I(s/4, v^{1/2}) = \int_{s/4}^{v^{1/2}} H_U(u,\mathcal{F})^{1/2}du \leq \frac{\theta^{3/2}n^{1/2}M}{256}. \tag{A.3}$$

*Then*

$$P^*(\sup_{h\in\mathcal{F}}|P_n h - Ph| > 4M) \leq 10\Big(1-\frac{1}{32\psi(M,n,v)}\Big)^{-1}\exp(-(1-\theta)\psi(M,n,v)).$$

**Proof:** The proof uses conditioning and chaining. The first step is conditioning. Let $\boldsymbol{Z}_1,\dots,\boldsymbol{Z}_N$ be an i.i.d. sample from $P$, and let $(R_1,\dots,R_N)$ be uniformly distributed over the set of permutations of $(1,\ldots,N)$, where $N = mn$, with $m = 2$. Define $n' = N - n$, $\tilde{P}_{n,N} = n^{-1}\sum_{i=1}^{n}\delta_{\boldsymbol{Z}_{R_i}}$, and $P_N = N^{-1}\sum_{i=1}^{N}\delta_{\boldsymbol{Z}_i}$, with $\delta_{\boldsymbol{Z}_i}$ the Dirac measure at observation $\boldsymbol{Z}_i$. Then the following inequality can be thought of as an alternative to the classical symmetrization inequality (cf., (25) Lemma 2.14.18 with $a = 2^{-1}$ and $m = 2$),

$$P^*(\sup_{h\in\mathcal{F}}|P_n h - Ph| > 4M) \le \Big(1 - \frac{v}{4nM^2}\Big)^{-1} P^*(\sup_{\mathcal{F}}|\tilde{P}_{n,N}h - P_N h| > M). \tag{A.4}$$

Conditioning on $\boldsymbol{Z}_1,\ldots,\boldsymbol{Z}_N$, it suffices to consider $P^*_{|N}(\sup_{\mathcal{F}}|\tilde{P}_{n,N}h - P_N h| > M)$, where $P_{|N}$ be the conditional distribution given $\boldsymbol{Z}_1,\ldots,\boldsymbol{Z}_N$.

The second step is to bound $P^*_{|N}(\sup_{\mathcal{F}}|\tilde{P}_{n,N}h - P_N h| > M)$ by chaining. Let $\varepsilon_0 > \varepsilon_1 > \ldots > \varepsilon_T > 0$ be a sequence to be specified. Denote by $\mathcal{F}_q$ the minimal $\varepsilon_q$-net for $\mathcal{F}$ with respect to the $L_2(P_N)$-norm. For each $h$, let $\pi_q h = \arg\min_{g\in\mathcal{F}_q}\|g - h\|_{P_N,2}$. Evidently, $\|\pi_q h - h\|_{P_N,2} \le \varepsilon_q$, and $|\mathcal{F}_q| = N(\varepsilon_q,\mathcal{F},L_2(P_N))$, the covering number. Then $P^*_{|N}(\sup_{\mathcal{F}}|\tilde{P}_{n,N}h - P_N h| > M)$ is bounded by

$$\begin{aligned}
& P^*_{|N}(\sup_{\mathcal{F}}|(\tilde{P}_{n,N} - P_N)(\pi_0 h)| > (1-\frac{\theta}{4})M) \\
& \qquad + P^*_{|N}(\sup_{\mathcal{H}}|(\tilde{P}_{n,N} - P_N)(\pi_0 h - \pi_T h)| > \frac{\theta M}{8}) \\
& \qquad + P^*_{|N}(\sup_{\mathcal{F}}|(\tilde{P}_{n,N} - P_N)(\pi_T h - h)| > \frac{\theta M}{8}) \\
\le\ & |\mathcal{F}_0|\sup_{\mathcal{F}} P^*_{|N}(|(\tilde{P}_{n,N} - P_N)(\pi_0 h)| > (1-\frac{\theta}{4})M) \\
& + \sum_{q=1}^{T}|\mathcal{F}_q||\mathcal{F}_{q-1}|\sup_{\mathcal{F}} P^*_{|N}(|(\tilde{P}_{n,N} - P_N)(\pi_q h - \pi_{q-1}h)| > \eta_q) \\
& + \sup_{\mathcal{F}} P^*_{|N}(|(\tilde{P}_{n,N} - P_N)(\pi_T h - h)| > \frac{\theta M}{8}) \\
:=\ & P_1 + P_2 + P_3,
\end{aligned}$$

where

$$\eta_q = \varepsilon_{q-1}(\frac{16H_U(\varepsilon_q,\mathcal{F})}{\theta n})^{1/2}; q = 1,\ldots,T, \tag{A.5}$$

and $\varepsilon_0 = H_U(\frac{\theta}{4}\psi(M,n,v),\mathcal{F})^-$, $\varepsilon_{q+1} = s \vee \sup\{x \le \varepsilon_q/2 : H_U(x,\mathcal{F}) \ge 4H_U(\varepsilon_q,\mathcal{F})\}$; $q = 0,\ldots,T$, and $T = \min\{q : \varepsilon_q \le s\}$. Note that $\varepsilon_0 \le v^{1/2}$ by construction. Furthermore, by (A.3) and Lemma 3.1 of (1),

$$\sum_{q=1}^{T}\eta_j = \sum_{q=1}^{T}\varepsilon_{q-1}(\frac{16H_U(\varepsilon_q,\mathcal{F})}{\theta n})^{1/2} \le \frac{32}{(\theta n)^{1/2}}I(s/4,v^{1/2}) \le \theta M/8. \tag{A.6}$$

We now proceed to bound $P_1$-$P_3$ separately.

On $C_N = (\sup_{\mathcal{F}} P_N h^2 \leq 64v)$, $\sigma_N^2 = P_N(\pi_0 f - P_N \pi_0)^2 \leq P_N(\pi_0 f)^2 \leq 64v$, by Massart's inequality, cf., (25), Lemma 2.14.19, $P^*_{|N}(|(\tilde{P}_{n,N} - P_N)(\pi_0 h)| > (1 - \frac{\theta}{4})M) \leq 2\exp(-n(1-\theta/4)^2 M^2/(2\sigma_N^2)) \leq 2\exp(-(1-\theta/4)^2\psi(M,n,v))$. By the choice of $\varepsilon_0$, $P_1 \leq 2\exp(H_U(\varepsilon_0, \mathcal{F}))\exp(-(1-\theta/4)^2\psi(M,n,v)) \leq 2\exp(-(1-\theta)\psi(M,n,v))$. On $C_N^c$, it follows from Lemma 33 of (18) that $P^*(C_N^c)$ is bounded by $P^*(\sup_{\mathcal{F}}(P_N h^2)^{1/2} \geq 8v^{1/2}) \leq 4\exp(-Nv + H_U(v^{1/2}, \mathcal{F})) \leq 4\exp(-2nv + \frac{\theta}{4}\psi(M,n,v)) \leq 4\exp(-(1-\theta)\psi(M,n,v))$

For $P_2$, if $\epsilon_0 \leq s$, let $\varepsilon_T = \varepsilon_0$. Then $P_2 = 0$. Otherwise, consider the case of $\epsilon_0 > s$. Note that $P_N(\pi_q h - \pi_{q-1}h)^2 \leq 2(P_N(\pi_q h - h)^2 + P(h - \pi_{q-1}h)^2) \leq 2\varepsilon_q^2 + 2\varepsilon_{q-1}^2 \leq 4\varepsilon_{q-1}^2$. By Massart's inequality, $P^*_{|N}(|(\tilde{P}_{n,N} - P_N)(\pi_q h - \pi_{q-1}h)| > \eta_q) \leq 2\exp(-n\eta_q^2/(2\sigma_N^2))$ with $\sigma_N^2 \leq P_N(\pi_q h - \pi_{q-1}h)^2 \leq 4\varepsilon_{q-1}^2$, and by the choice of $\eta_q$, $q = 1, \ldots, T$,

$$
\begin{aligned}
P_2 &\leq \sum_{q=1}^{N} |\mathcal{F}_q|^2 \sup_{\mathcal{F}} P^*_{|N}(|(\tilde{P}_{n,N} - P_N)(\pi_q h - \pi_{q-1}h)| > \eta_q) \\
&\leq 2\sum_{q=1}^{T} \exp(2H_U(\varepsilon_q, \mathcal{F}) - \frac{n\eta_q^2}{4m\varepsilon_{q-1}^2}) = 2\sum_{q=1}^{T} \exp((2 - 2/\theta)H_U(\varepsilon_q, \mathcal{F})) \\
&\leq 2\sum_{q=1}^{\infty} \exp((2 - 2/\theta)4^q H_U(\varepsilon_0, \mathcal{F})) \leq 4\exp(-(1-\varepsilon)\psi(M,n,v)).
\end{aligned}
$$

For $P_3$, note that $\tilde{P}_{n,N} f \leq 2P_N f$ for any $f \geq 0$, and $P_N(\pi_T h - h)^2 \leq \varepsilon_T^2$ by the definition of $\pi_T$. Then $|(\tilde{P}_{n,N} - P_N)(\pi_T h - h)|^2 \leq 2(\tilde{P}_{n,N} + P_N)(\pi_T h - h)^2 \leq 6\varepsilon_T^2 \leq (\theta M/8)^2$ because $\varepsilon_T \leq s = \frac{\theta M}{8\sqrt{6}}$. So $P_3 = 0$.

Now

$$
\begin{aligned}
P|N^*(\sup_{\mathcal{F}} |\tilde{P}_{n,N} h - P_N h| > M) &\leq P_1 + P_2 + P_3 \\
&\leq 6\exp(-(1-\theta)\psi(M,n,v)).
\end{aligned}
$$

After taking the expectation with respect to $\boldsymbol{Z}_1, \ldots, \boldsymbol{Z}_N$, we have, from (A.4), that $P^*(\sup_{\mathcal{F}} |P_n h - Ph| > 4M)$ is upper bounded by

$$
10\Big(1 - \frac{1}{32\psi(M,n,v)}\Big)^{-1} \exp(-(1-\theta)\psi(M,n,v)).
$$

This completes the proof.

**Proof of Lemma 2 :** It can be verified that, with $\lambda_j > 0$ constants,

$$
\begin{aligned}
R_{V_1}(a,b) &= \lambda_1 \left(\theta_1 e^{-b} \int_0^1 e^{-ax} x^\gamma dx + (1-\theta_1)e^b \int_0^1 e^{ax} x^\gamma dx + \right. \\
&\quad \left. \theta_2 e^{-b} \int_{-1}^0 e^{-ax}(-x)^\gamma dx + (1-\theta_2)e^b \int_{-1}^0 e^{ax}(-x)^\gamma dx \right);
\end{aligned}
$$

$R_{V_2}(a,b)$ can be expressed as

$$\lambda_2\Big(\theta_1\int_0^1\log(1+e^{-ax-b})x^\gamma dx+(1-\theta_1)\int_0^1\log(1+e^{ax+b})x^\gamma dx+$$
$$\theta_2\int_{-1}^0\log(1+e^{-ax-b})(-x)^\gamma dx+(1-\theta_2)\int_{-1}^0\log(1+e^{ax+b})(-x)^\gamma dx\Big);$$

$R_{V_3}(a,b)$ can be written, in the region of interest $\{(a,b): a\le 0, -1\le -(1+b)a^{-1}\le 0, 0\le (1-b)a^{-1}\le 1\}$, as

$$\lambda_3(\theta_1((\gamma+1)(\gamma+2))^{-1}(1-b)^{\gamma+2}a^{-\gamma-1}+(1-\theta_1)(\frac{1+b}{\gamma+1}+\frac{a}{\gamma+2})+$$
$$\theta_2(\frac{1-b}{\gamma+1}+\frac{a}{\gamma+2})+(1-\theta_2)((\gamma+1)(\gamma+2))^{-1}(1+b)^{\gamma+2}a^{-\gamma-1});$$

$R_{V_4}(a,b)$ can be written as, when $a>0$ and $b>1$,

$$\lambda_4\Big(\frac{\theta_1}{(\gamma+1)(\gamma+2)}(\frac{(1-b)^{\gamma+2}}{a^{\gamma+1}}-\frac{(-b)^{\gamma+2}}{a^{\gamma+1}})+(1-\theta_1)$$
$$(\frac{1}{\gamma+1}-\frac{1}{(\gamma+1)(\gamma+2)}(\frac{(-b)^{\gamma+2}}{a^{\gamma+1}}-\frac{(-b-1)^{\gamma+2}}{a^{\gamma+1}}))+\theta_2/(\gamma+1)\Big);$$

when $a>0$ and $0<b\le 1$,

$$\lambda_4\Big(\frac{\theta_1}{(\gamma+1)(\gamma+2)}(\frac{(1-b)^{\gamma+2}}{a^{\gamma+1}}-\frac{(-b)^{\gamma+2}}{a^{\gamma+1}})+(1-\theta_1)(\frac{1}{\gamma+1}$$
$$-\frac{1}{(\gamma+1)(\gamma+2)}(\frac{(-b)^{\gamma+2}}{a^{\gamma+1}}))+\theta_2/(\gamma+1)+\frac{1-\theta_2}{(\gamma+1)(\gamma+2)}\frac{(1+b)^{\gamma+2}}{a^{\gamma+1}}\Big);$$

when $a>0$ and $-1<b\le 0$,

$$\lambda_4\Big(\frac{\theta_1}{(\gamma+1)(\gamma+2)}\frac{(1-b)^{\gamma+2}}{a^{\gamma+1}}+(1-\theta_1)/(\gamma+1)+\theta_2(\frac{1}{\gamma+1}$$
$$-\frac{1}{(\gamma+1)(\gamma+2)}(\frac{b^{\gamma+2}}{a^{\gamma+1}}))+\frac{1-\theta_2}{(\gamma+1)(\gamma+2)}(\frac{(1+b)^{\gamma+2}}{a^{\gamma+1}}-\frac{b^{\gamma+2}}{a^{\gamma+1}})\Big);$$

when $a>0$ and $b\le -1$,

$$(1-\theta_1)/(\gamma+1)+\theta_2(\frac{1}{\gamma+1}-\frac{1}{(\gamma+1)(\gamma+2)}(\frac{b^{\gamma+2}}{a^{\gamma+1}}-\frac{(b-1)^{\gamma+2}}{a^{\gamma+1}}))$$
$$+\frac{1-\theta_2}{(\gamma+1)(\gamma+2)}(\frac{(1+b)^{\gamma+2}}{a^{\gamma+1}}-\frac{b^{\gamma+2}}{a^{\gamma+1}})).$$

Similarly $R(a,b)=\frac{1}{2}((1-\theta_1+\theta_2)+(I[b<0](2\theta_1-1)+I[b>0](1-2\theta_2))|b/a|^{\gamma+1})$ when $a>0$. The results can be verified through direct calculation.

**Proof of Lemma 3 :** For $V_j$, $j=1,\cdots,3$, let $\boldsymbol{f}^*=f^{V_j}$. We verify Assumptions A and B through the exact expression of $V_j$ given in the proof of Lemma 2, although Taylor's expansion is generally applicable. Note $R_{V_j}(a,b)$ is strictly convex and smooth; $j=1,2$ and $R_{V_3}(a,b)$ is piecewise smooth and strictly convex in the neighborhood of $(a_3^*,b_3^*)$. For any $(a,b)$ in the neighborhood of $(a_j^*,b_j^*)$ and some constant $d_1>0$, $e_{V_j}(\boldsymbol{f},\boldsymbol{f}^{V_j}) = R_{V_j}(a,b) - R_{V_j}(a_j^*,b_j^*) \geq d_1(a-a_j^*,b-b_j^*)H_{V_j}(a_j^*,b_j^*)(a-a_j^*,b-b_j^*)^T$ with a positive definite matrix $H_{V_j}(a_j^*,b_j^*)$. Moreover, $e(\boldsymbol{f},\boldsymbol{f}^{V_j}) = R(a,b)-R(a_j^*,b_j^*) = (I[b<0](2\theta_1-1)+I[b>0](1-2\theta_2))(|b/a|^{\gamma+1}-|b_j^*/a_j^*|^{\gamma+1})/2$. By the assumption that $\theta_1+\theta_2\neq 1$, $b_j^*\neq 0$, $|e(\boldsymbol{f},\boldsymbol{f}^{V_j})| \leq c_1|b_j/a_j - b_j^*/a_j^*|$ for some constant $c_1>0$, implying Assumption A with $\alpha=1/2$. For Assumption B, note that $|f|\leq T_1$ for some constant $T_1>0$ when $e_V(\boldsymbol{f},\boldsymbol{f}^{V_j})$ is small. Then $|V_j(\boldsymbol{f},\boldsymbol{z})-V_j(\boldsymbol{f}^{V_j},\boldsymbol{z})| \leq V_j'(-T_1)|f(x)-f^{V_j}(x)| = V_j'(-T_1)|(a-a_j^*)x+(b-b_j^*)|$. Hence $Var(V_j(\boldsymbol{f},Z)-V_j(\boldsymbol{f}^{V_j},Z)) \leq c_2E((a-a_j^*)X+(b-b_j^*))^2 = (a-a_j^*,b-b_j^*)D_{V_j}(a-a_j^*,b-b_j^*)^T$ with $D_{V_j}$ a positive definite matrix, implying Assumption B with $\beta=1$.

For $V_4$, the minimal of $R_{V_4}(a,b)$ is attained as $a\to\infty$ and $b=0$, independent of $\theta_1$, $\theta_2$ and $\gamma$. Note that $\boldsymbol{f}^{V_4}=\bar{\boldsymbol{f}}$ since $EV_4(\bar{\boldsymbol{f}},Z)=\inf_{\boldsymbol{f}\in\mathcal{F}} EV_4(\boldsymbol{f},Z)$. Direct calculation yields that $e_V(\boldsymbol{f},\bar{\boldsymbol{f}}) = R_{V_4}(a,b) - \lim_{a\to\infty} R_{V_4}(a,0) \geq c_4|b/a|^{\gamma+1}$, and $|e(\boldsymbol{f},\bar{\boldsymbol{f}})| \leq c_5|b/a|^{\gamma+1}$. This implies Assumption A with $\alpha=1$. For Assumption B, it follows from the fact that $V_4(\bar{\boldsymbol{f}},\boldsymbol{z}) = L(\bar{\boldsymbol{f}},\boldsymbol{z})$ for any $\boldsymbol{z}$ that $Var(V_4(\boldsymbol{f},Z)-V_4(\bar{\boldsymbol{f}},Z)) \leq 2E|V_4(\boldsymbol{f},Z)-V_4(\bar{\boldsymbol{f}},Z)| \leq 2E|L(\boldsymbol{f},Z)-L(\bar{\boldsymbol{f}},Z)| + 2E(V_4(\boldsymbol{f},Z)-L(\boldsymbol{f},Z))$. Furthermore, $E(L(\boldsymbol{f},Z)-L(\bar{\boldsymbol{f}},Z)) = E|2P(Y=1|X=x)-1||L(\boldsymbol{f},Z)-L(\bar{\boldsymbol{f}},Z)| \geq \min(|2\theta_1-1|,|2\theta_2-1|)E|L(\boldsymbol{f},Z)-L(\bar{\boldsymbol{f}},Z)|$, and $E(V_4(\boldsymbol{f},Z)-L(\boldsymbol{f},Z)) \leq E(V_4(\boldsymbol{f},Z)-L(\bar{\boldsymbol{f}},Z)) = e_{V_4}(\boldsymbol{f},\bar{\boldsymbol{f}})$. Therefore, Assumption B is met with $\beta=1$. This completes the proof.

**Proof of Lemma 4:** For Assumption A, note that for any $\boldsymbol{w}=\boldsymbol{w}^*+\Delta\boldsymbol{w}$ in a small neighborhood of $\boldsymbol{w}^*$, $e_V(\boldsymbol{f},\boldsymbol{f}^V) = R_V(\boldsymbol{w})-R_V(\boldsymbol{w}^*) \geq d_2\Delta\boldsymbol{w}^T H_1\Delta\boldsymbol{w}$. Furthermore, direct computation yields that $e(\boldsymbol{f},\boldsymbol{f}^V) \leq d_3(\Delta\boldsymbol{w}^T\Delta\boldsymbol{w})^{(\gamma+1)/2}$, for some constant $d_3>0$. Hence, for some constant $c_2>0$, $|e(\boldsymbol{f},\boldsymbol{f}^V)| \leq c_2e_V(\boldsymbol{f},\boldsymbol{f}^V)^{(\gamma+1)/2}$ for all small $\Delta\boldsymbol{w}$, implying Assumption A with $\alpha=(\gamma+1)/2$. For Assumption B, it follows from the fact $|V(\boldsymbol{f},\boldsymbol{z})-V(\boldsymbol{f}^V,\boldsymbol{z})| \leq \sum_{c=1}^4|f_c(x)-f_c^V(x)|$ that $Var(V(\boldsymbol{f},Z)-V(\boldsymbol{f}^V,Z))$ is upper bounded by

$$E\Big(\sum_{c=1}^4|f_c(X)-f_c^V(X)|\Big)^2 \leq 4E\Big(\sum_{c=1}^4(f_c(X)-f_c^V(X))^2\Big) = \Delta\boldsymbol{w}^T H_2\Delta\boldsymbol{w}$$

with $H_2$ a positive definite matrix, implying that $Var(V(\boldsymbol{f},Z)-V(\boldsymbol{f}^V,Z)) \leq c_3e_V(\boldsymbol{f},\boldsymbol{f}^V)$ for some constant $c_3>0$ and all small $\Delta\boldsymbol{w}$, and thus Assumption B with $\beta=1$. This completes the proof.

**Proof of Lemma 5:** We use a pointwise argument. First consider $V=V_{svm1}$. Note that $EV(\boldsymbol{f},\boldsymbol{Z}) = E(E(\sum_{c=1}^3 p_c(X)h_{svm1}(\boldsymbol{u}(\boldsymbol{f}(X),c))|X))$. Now define $h_{\boldsymbol{p}}(\boldsymbol{f})$ to be $\sum_{c=1}^3 p_cV_c(\boldsymbol{f})$ for any $\boldsymbol{f}\in\mathbb{R}^3$, where $V_c(\boldsymbol{f}) = \sum_{c=1}^3 p_c\sum_{j\neq c}(1-(f_c-f_j))_+$. We now verify that $\bar{\boldsymbol{f}}=(2/3,-1/3,-1/3)$ minimizes $h_{\boldsymbol{p}}(\boldsymbol{f})$ when $\boldsymbol{p}=(5/11,3/11,3/11)$, for $x\in[0,1/3]$. The other two cases when $x\in(1/3,2/3]$ and $x\in(2/3,1]$ can be dealt with similarly. Now reparametrize $\boldsymbol{f}$ as $\bar{\boldsymbol{f}}+(r_1,r_2,r_3)^Td$ with $\sum_c r_c=0$, $\sum_c|r_c|=1$ and $d=\|\boldsymbol{f}-\bar{\boldsymbol{f}}\|_1$. When $d$ is suffi-ciently small and $\boldsymbol{p} = (5/11, 3/11, 3/11)$,

$$\begin{aligned} h_{\boldsymbol{p}}(\boldsymbol{f}) - h_{\boldsymbol{p}}(\bar{\boldsymbol{f}}) &= p_1((r_2 - r_1)_+ d + (r_3 - r_1)_+ d) + p_2(2 + (r_1 - r_2)d + 1 \\ &\quad +(r_3 - r_2)d) + p_3(2 + (r_1 - r_3)d + 1 + (r_2 - r_3)d) \\ &\quad -(3p_2 + 3p_3) \\ &= (p_1(r_2 - r_1)_+ + p_1(r_3 - r_1)_+ + p_2(r_1 - r_2) + p_2(r_1 - r_3))d \\ &\geq C_1 d, \end{aligned}$$

where $C_1 = \min_{\sum_c r_c = 0, \sum_c |r_c| = 1} (p_1(r_2 - r_1)_+ + p_1(r_3 - r_1)_+ + p_2(r_1 - r_2) + p_2(r_1 - r_3)) > 0$. By convexity of $h_{\boldsymbol{p}}(\boldsymbol{f})$, $\bar{\boldsymbol{f}}$ is the minimizer. Combining the three cases, we obtain with $\bar{\boldsymbol{f}}(x) = (I[0 < x \leq 1/3] - 1/3, I[1/3 < x \leq 2/3] - 1/3, I[2/3 \leq x \leq 1] - 1/3)$ that $\boldsymbol{f}^*(x)$ minimizes $h_{p(x)}(\boldsymbol{f}(x))$ for each $x$, thus $EV(\boldsymbol{f}, Z) = Eh_{p(X)}(\boldsymbol{f}(X))$. For $V = V_\psi$, it follows from an argument similar to the proof of Theorem 3.2 of (16). This completes the proof.

**Proof of Lemma 6:** We will apply an argument similar to that in the proof of Lemma 5. Let $h^T_{\boldsymbol{p}}(\boldsymbol{f})$ be $\sum_{c=1}^3 p_c V^T_c(\boldsymbol{f})$ with $V^T(\boldsymbol{f}) = T \wedge V_c(\boldsymbol{f})$ with $V_c(\boldsymbol{f})$ defined in the proof of Lemma 5. Note that $f^* = \bar{\boldsymbol{f}}$ and $e_{V^T}(\boldsymbol{f}, \boldsymbol{f}^*) = E(E(V^T(\boldsymbol{f}, \boldsymbol{Z}) - V(\boldsymbol{f}^*, \boldsymbol{Z})|X)) = E(h_{\boldsymbol{p}(X)}(\boldsymbol{f}(X)) - h_{\boldsymbol{p}(X)}(\boldsymbol{f}^*(X)))$. By Theorem 3.1 of Liu and Shen (2006), $e(\boldsymbol{f}, \boldsymbol{f}^*) = E(\max_c p_c(X) - p_{\text{argmin}_c f_c(X)}(X))$. It then suffices to show that $h_{\boldsymbol{p}}(\boldsymbol{f}) - h_{\boldsymbol{p}}(\boldsymbol{f}^*) \geq C(\max_c p_c - p_{\text{argmin}_c f_c})$ for any measurable $\boldsymbol{f} \in \mathbb{R}^3$ and some constant $C > 0$. Suppose $\boldsymbol{p} = (5/11, 3/11, 3/11)$ without loss of generality. Then $p_1 = \max_c p_c$, and the proof becomes trivial when $f_1 = \max_c f_c$. Suppose $f_1 < \max_c f_c$ and further $f_1 < f_2$ without loss of generality. Then $\|\boldsymbol{f} - \boldsymbol{f}^*\|_1 \geq |f_1 - f_1^*| + |f_2 - f_2^*| \geq |(f_1^* - f_2^*) + (f_2 - f_1)| > 1$ ($f_1^* = 2/3$ and $f_2^* = -1/3$). Two cases are treated separately. When $\max_c V_c(\boldsymbol{f}) < T$ or $V_c^T(\boldsymbol{f}) = V_c(\boldsymbol{f})$, $h^T_{\boldsymbol{p}}(\boldsymbol{f}) - h^T_{\boldsymbol{p}}(\boldsymbol{f}^*) \geq C_1 \|\boldsymbol{f} - \boldsymbol{f}^*\|_1 \geq C_1$ by the proof of Lemma 5. When $\max_c V_c(\boldsymbol{f}) \geq T$ or $\max_c V_c^T(\boldsymbol{f}) = T$, $h^T_{\boldsymbol{p}}(\boldsymbol{f}) \geq (\min_c p_c)(\sum_{c=1}^3 V_c^T(\boldsymbol{f})) \geq (3/11)T$, and $h_{\boldsymbol{p}}(\boldsymbol{f}^*) = 18/11$. Then, $h_{\boldsymbol{p}}(\boldsymbol{f}) - h_{\boldsymbol{p}}(\boldsymbol{f}^*) \geq 3T/11 - 18/11 = 9/11 > \max_c p_c - p_{\text{argmin}_c f_c} = 2/11$. The desired result follows.

**Proof of Lemma 7:** The proof uses the pointwise argument and is similar to that of Lemma 6. Note that $E(V^T(\boldsymbol{f}, \boldsymbol{Z}) - V(\boldsymbol{f}^*, \boldsymbol{Z}))^2 \leq TE|V^T(\boldsymbol{f}, \boldsymbol{Z}) - V(\boldsymbol{f}^*, \boldsymbol{Z})| = TE(p_1(X)|V_1^T(\boldsymbol{f}(X)) - V_1(\boldsymbol{f}(X))| + p_2(X)|V_2^T(\boldsymbol{f}(X)) - V_2(\boldsymbol{f}^*(X))| + p_3(X)|V_3^T(\boldsymbol{f}(X)) - V_3(\boldsymbol{f}^*(X))|)$, with $V_c(\boldsymbol{f})$ defined in the proof of Lemma 5. It suffices to show that for any measurable $\boldsymbol{f} \in \mathbb{R}^3$ $Left \leq CRight$ for some constant $C > 0$ with $Left \equiv p_1|V_1^T(\boldsymbol{f}) - V_1(\boldsymbol{f}^*)| + p_2|V_2^T(\boldsymbol{f}) - V_2(\boldsymbol{f}^*)| + p_3|V_3^T(\boldsymbol{f}) - V_3(\boldsymbol{f}^*)|$ and $Right \equiv p_1(V_1^T(\boldsymbol{f}) - V_1(\boldsymbol{f}^*)) + p_2(V_2^T(\boldsymbol{f}) - V_2(\boldsymbol{f}^*)) + p_3(V_3^T(\boldsymbol{f}) - V_3(\boldsymbol{f}^*))$. Two cases are examined.

(1) If $M(x) \equiv \max_{1 \leq c \leq 3} V_c(\boldsymbol{f}) \geq T$, then $\max_{1 \leq c \leq 3} V_c^T(\boldsymbol{f}) = T$. It follows that $Left \leq T$ and $Right \geq 3T/11 - 18/11$, as shown in the proof of Lemma 6. This implies that $Left \leq 11 Right$ because $T \geq 9$.

(2) If $M(x) \leq T$, we prove the non-truncated version of the inequality. Note that $|V_c(\boldsymbol{f}) - V_c(\boldsymbol{f}^*)| \leq \|\boldsymbol{f} - \boldsymbol{f}^*\|_1$, $c = 1, 2, 3$. Then $Left \leq \|\boldsymbol{f} - \boldsymbol{f}^*\|_1$ for any $\boldsymbol{f} \in \mathbb{R}^3$. Following the proof of Lemma 6, we have $Right > C_1 \|\boldsymbol{f} - \boldsymbol{f}^*\|_1 \geq C_1 Left$.

**Proof of Lemma 8:** Note that $H_U(2s\varepsilon, \mathcal{G}_p(s)) = H_U(\varepsilon, \mathcal{G}_p^*(s))$ with $\mathcal{G}_p^*(s) = \{(2s)^{-1}f : f \in \mathcal{G}_p(s)\}$. In addition, $\mathcal{G}_p^*(s)$ is the convex hull of $\pm x_j/2$, $j = 1, \dots, p$. Then $H_U(k^{-1/2}, \mathcal{G}_p^*(s)) \leq \log\binom{2p+k-1}{k} \leq \log((2p+k)!) - \log((2p)!) - \log(k!)$ for any integer $k > 1$ using the argument in the proof of Lemma 2.6.11 of (25). By Stirling's formula,

$$\sqrt{2\pi}n^{n+1/2}e^{-n+1/(12n+1)} < n! < \sqrt{2\pi}n^{n+1/2}e^{-n+1/(12n)},$$

implying that $H_U(k^{-1/2}, \mathcal{G}_p^*(s))$ is no greater than $2p\log(1 + k/(2p))$ $+ k\log(2p/k + 1) - \log\sqrt{2\pi} + (12(2p+k))^{-1} - (12(2p)+1)^{-1} - (12k+1)^{-1}$. Let $\varepsilon > k^{-1/2}$. Then there exists a constant $c > 0$ such that $H_U(\varepsilon, \mathcal{G}_p^*(s)) \leq c(p\log(1 + \frac{1}{p\varepsilon^2}) + \varepsilon^{-2}\log(p\varepsilon^2 + 1))$. This completes the proof.